\documentclass[runningheads]{llncs}
\usepackage[T1]{fontenc}
\usepackage{stmaryrd}
\usepackage{mathptmx}

\usepackage{float}

\usepackage{xcolor}
\usepackage{colortbl}

\usepackage{amsmath,amssymb}

\spnewtheorem{thm}{Theorem}[section]{\bfseries}{\itshape}
\spnewtheorem{lem}[thm]{Lemma}{\bfseries}{\itshape}
\spnewtheorem{cor}[thm]{Corollary}{\bfseries}{\itshape}
\spnewtheorem{defn}[thm]{Definition}{\bfseries}{\itshape}
\spnewtheorem{prop}[thm]{Proposition}{\bfseries}{\itshape}
\spnewtheorem{fact}[thm]{Fact}{\bfseries}{\itshape}
\spnewtheorem{exm}[thm]{Example}{\bfseries}{\itshape}



\usepackage{amscd}

\usepackage{multirow}

\DeclareFontFamily{U}  {MnSymbolC}{}
  
  \DeclareFontShape{U}{MnSymbolC}{m}{n}{
    <-6>  MnSymbolC5
   <6-7>  MnSymbolC6
   <7-8>  MnSymbolC7
   <8-9>  MnSymbolC8
   <9-10> MnSymbolC9
  <10-12> MnSymbolC10
  <12->   MnSymbolC12}{}
\DeclareFontShape{U}{MnSymbolC}{b}{n}{
    <-6>  MnSymbolC-Bold5
   <6-7>  MnSymbolC-Bold6
   <7-8>  MnSymbolC-Bold7
   <8-9>  MnSymbolC-Bold8
   <9-10> MnSymbolC-Bold9
  <10-12> MnSymbolC-Bold10
  <12->   MnSymbolC-Bold12}{}

\DeclareSymbolFont{MnSyC} {U} {MnSymbolC}{m}{n}
\DeclareMathSymbol{\veedot}{\mathop}{MnSyC}{47}

\usepackage{rotating}

\usepackage{relsize}

\usepackage{graphicx}
\usepackage{color}
\usepackage{xspace}
\usepackage{cleveref}

\crefname{thm}{theorem}{theorems}
\Crefname{thm}{Theorem}{Theorems}
\crefname{lem}{lemma}{lemmas}
\Crefname{lem}{Lemma}{Lemmas}
\crefname{prop}{proposition}{propositions}
\Crefname{prop}{Proposition}{Propositions}
\crefname{cor}{corollary}{corollaries}
\Crefname{cor}{Corollary}{Corollaries}
\crefname{defn}{definition}{definitions}
\Crefname{defn}{Definition}{Definitions}
\crefname{fact}{fact}{facts}
\Crefname{fact}{Fact}{Facts}
\crefname{exm}{example}{examples}
\Crefname{exm}{Example}{Examples}

\usepackage{bussproofs}

\usepackage{tikz}
\usepackage{bpextra}

\usepackage{enumerate}
\usepackage{array}
\newcolumntype{C}[1]{>{\centering\arraybackslash}p{#1}}

\DeclareSymbolFont{letters}{OML}{cmm}{m}{it}

\DeclareMathAlphabet{\mathcal}{OMS}{cmsy}{m}{n}

\newcommand{\PL}{\ensuremath{\mathbf{PL}}\xspace}

\newcommand{\LL}{\ensuremath{\mathsf{L}}\xspace}

\newcommand{\Prop}{\ensuremath{\mathsf{Prop}}\xspace}

\newcommand{\dblsetminus}{\mathbin{{\setminus}\mspace{-5mu}{\setminus}}}

\newcommand{\dep}{\ensuremath{\mathop{=\!}}\xspace}

\newcommand{\vvee}{\raisebox{1pt}{\ensuremath{\,\mathop{\mathsmaller{\mathsmaller{\dblsetminus\hspace{-0.23ex}/}}}\,}}}
\newcommand{\bigvvee}{\ensuremath{\mathop{\mathlarger{\dblsetminus}\hspace{-0.28ex}\raisebox{0.3pt}{$\mathlarger{{/}}$}}}\xspace}

\newcommand{\dneg}{\ensuremath{\mathop{-}}\xspace}

\newcommand{\barrightharpoon}{%
    \mathrel{\vcenter{\offinterlineskip
    \halign{\hfil##\hfil\cr
    \(\overline{\kern-0.1em\rightharpoondown\kern-0.1em}\)\cr
    \cr}}}}

\begin{document}
\title{%Dependence and inclusion over propositional teams
There are (other) ways to negate in propositional team semantics%: a note
}
\author{Fan Yang%\inst{1}%\orcidID{0000-0003-0392-6522} %\and
%Second Author\inst{2,3}\orcidID{1111-2222-3333-4444} \and
%Third Author\inst{3}\orcidID{2222--3333-4444-5555}
}
\authorrunning{F. Yang}
\institute{Department of Philosophy and Religious Studies, Utrecht University,
%Janskerkhof 13, 3512 BL Utrecht, 
The Netherlands
\email{f.yang@uu.nl}
}
\maketitle              % typeset the header of the contribution
\begin{abstract}
%The abstract should briefly summarize the contents of the paper in
%Negation is often a tricky connective for logics based on team semantics. The syntax of these logics thus usually only allows atomic negation or restricted negation. 
The languages of logics based on team semantics typically only allow  atomic negation or restricted negation.
In this paper, we explore propositional team-based logics with full (intuitionistic) negation. We demonstrate that including full intutionistic negation does not complicate the axiomatization of propositional  team-based logics with the downward closure property. We also review  known expressive completeness results for these logics, highlighting how relevant complemented properties are expressed in propositional dependence logic without directly using negation. Building on these insights, we also prove a new result: propositional logic extended with both dependence and inclusion atoms is expressively complete.

%We revisit the proofs of expressive completeness results for several propositional team-based logics,  highlighting  alternative ways to express complements of properties without using negation. We also prove a new result that team-based propositional logic with dependence and inclusion atoms is expressively complete. We illustrate that including full intutionistic negation does not complicate the proof theory of several propositional downward closed team logics. We also document a folklore simpler system for $\PL(\vvee)$.

\end{abstract}
\section{Introduction}

Hodges \cite{Hodges1997a,Hodges1997b} introduced {\em team semantics} as a compositional semantics for Hintikka and Sandu's independence-friendly logic \cite{HinSan96,Hintikka98book}. Team semantics was then further developed in V\"{a}\"{a}n\"{a}nen's dependence logic \cite{Van07dl} with atoms of the type $\dep(x,y)$ expressing that the value of $y$ is determined by the value of $x$. Team semantics characterises such dependency formulas by evaluating them over sets of assignments (called {\em teams}). This framework has since then been extended significantly. Many  logics based on team semantics have been introduced, such as independence logic (with independence atoms) \cite{D_Ind_GV}, inclusion logic (with inclusion atoms) \cite{Pietro_I/E,inclusion_logic_GH} and their propositional and modal variants. The idea of team semantics has also been independently adopted in inquisitive logic \cite{InquiLog} %which contains the inquisitive disjunction 
to characterise questions in natural language.

Many of these team-based logics do not include full negation in their syntax. Indeed, negation is often considered  a tricky connective in team semantics. The most commonly known negation in team semantics, the dual negation, originates from the game theoretic interpretation of team semantics. In the evaluation game of a formula, when a position associated with a negated formula $\dneg\alpha$ is reached, the Falsifier and Verifier swap their roles. %As observed already by Hodges \cite{Hodges1997a,Hodges1997b}, \todo{double check} 
 Such a negation is, however, not a semantic operator, in the sense that preservation of equivalence under replacement fails 
 %the meaning of equivalent formulas is not preserved 
 under the scope of negation. For instance, %in the context of first-order dependence logic, 
 we have $\dneg \dep(x,y)\equiv \bot$, whereas $\dneg\dneg \dep(x,y)\equiv \dep(x,y)\not\equiv\top\equiv\dneg\bot$. In fact, %dual negation's failure of being a semantic operator 
 the failure of dual negation to function as a semantic operator is even more extreme, as shown first by Burgess \cite{Burgess_negation_03} in the related context of Henkin quantifiers \cite{henkin61} and later by other authors in  %first-order, propositional and modal 
 team-based logics \cite{Negation_D_KV,Anttila2024}. This problem can, nevertheless, be avoided by considering both the team that validates the formula and the team that falsifies the formula (known also as ``co-team''), 
as done in Hodges' \cite{Hodges1997a,Hodges1997b} and in V\"{a}\"{a}n\"{a}nen's monograph \cite{Van07dl}. %, or alternatively, by using a bilateral approach to consider both \todo{rephrase, add }\cite{Aloni2022}

%atomic negation in team semantics is a type of dual negation coming from the original 

%The dual negation $\dneg$ that team-based logics inherited from the usual (single valuation-based) classical logic captures the role change of the two players in the semantic game of a formula. The semantics of this negation is not compositional, as observed already by Hodges \cite{Hodges1997a,Hodges1997b} in the context of IF-logic. For instance, in the context of first-order dependence logic, we have $\dneg \dep(p,q)\equiv \bot$, whereas $\dneg\dneg \dep(p,q)\equiv \dep(p,q)\not\equiv\dneg\bot\equiv \bot$. In fact, the failure of compositionality of the dual negation is even more dramatic, as shown by Burgess \cite{Burgess_negation_03} in the context of Henkin quantifiers and by other authors in the context of first-order, propositional and modal dependence logics \cite{Negation_D_KV,Anttila2024}. This problem can be fixed by considering team and co-team pairs in evaluating formulas, as done by Hodges in \cite{Hodges1997a,Hodges1997b} and in the monograph of V\"{a}\"{a}n\"{a}nen \cite{Van07dl}. In a similar manner, recent work by Aloni adopts a bilateral approach to include the dual or bilateral negation in team logics, whose semantics is defined by simultaneous induction for support and anti-support clauses.

Apart from these technical concerns, there is also a conceptual challenge with the dual negation. When defining the semantics of a negated dependence atom $\dneg\dep(x,y)$, in order to preserve the desirable downward closure property that all the other formulas in dependence logic have, one has to define $\dneg\dep(x,y)$ to be true only on the empty team (as $\dep(x,y)$ is always trivially true on any singleton team). In other words, $\dneg \dep(x,y)$ is defined to be equivalent to the constant falsum $\bot$. %This technical solution has often been criticised, 
This technical solution makes it difficult to provide a meaningful interpretation of the negated dependence atoms, and thus has often been criticised.

%the only  way to define a negated dependence atom is to set any such negation $\neg\dep(\vec{p},q)$ as true only on the empty team; in other words, $\neg \dep(\vec{p},q)$ is defined to be equivalent to the constant $\bot$, as $\dep(\vec{p},q)$ is always true on any singleton team. It is thus difficult to provide a meaningful interpretation for $\neg\dep(\vec{p},q)$. 

In order to avoid all these complications, most of the literature on team-based logics restricts the occurrence of  negation  in front of arbitrary formulas. This is done by assuming negation normal form and excluding negated dependence atoms from the class of well-formed formulas. 
However, negation as an important connective cannot be entirely avoided, particularly in the context of proof theory. In the partially complete deduction systems for first-order dependence, independence and inclusion logics \cite{Axiom_fo_d_KV,Hannula_fo_ind_13,YangInc20}\footnote{Note that these logics cannot be fully axiomatized, due to their strong expressive power (see the cited papers for details).}, the scope of negation was extended to cover all classical first-order formulas (i.e., those without dependence, independence or inclusion atoms). Similarly,  in propositional and modal team-based logics, negation is often allowed to occur in front of classical formulas (see, e.g., \cite{Yang2022,Yang17MD}). However, full negation is still excluded in these frameworks. In this paper, we study team-based propositional logics with full negation.

Negation can, of course, also be defined in terms of an appropriate implication. 
%One can of course define negation in terms of (an appropriate) implication. 
In particular, in inquisitive logic, the formula $\phi\to\bot$ (with the intuitionistic implication $\to$) is (as usual) treated as the definition of the (intuitionistic) negation $\neg\phi$ of $\phi$. In this paper, we define the full negation $\neg\phi$ exactly through the semantic clause of the intuitionistic negation $\phi\to\bot$.
%and this is essentially also the definition of the full negation we adopt in this paper. 
 This negation is conservative over the dual negation, as for formulas $\alpha$ in the original syntax of propositional logic, %classical formulas (formulas without inquisitive disjunction ), 
 $\neg\alpha$ (i.e., $\alpha\to\bot$) is equivalent to the dual negation $\dneg\alpha$. %We will present the preliminaries related to this full (intuitionistic) negation in Section 2. 
 %Let us already clarify that 
 Since such an intuitionistic negation 
  %On the other hand, failure of replacement does not occur for this intuitionistic negation, because its semantics 
 is  defined compositionally in the usual manner, without referring to co-teams, failure of replacement does not any more occur; %anti-support clauses, 
%for example). 
in particular, we have $\neg\dep(\vec{x},y)\equiv \bot$ and $\neg\neg\dep(\vec{x},y)\equiv \neg\bot\equiv \top$.

 In the context of dependence logic, while the intuitionistic implication (and thus intuitionsitic negation) have also been introduced  \cite{AbVan09}, they are usually not included in the standard syntax of first-order dependence logic and its variants. This is because, for instance, first-order dependence logic is expressively equivalent to existential second-order logic \cite{KontVan09,Van07dl}, whereas its extension with intuitionistic implication is expressively too strong---equivalent to full second-order logic \cite{Yang2011}--- and consequently many desirable properties, such as compactness and others, are lost.  In fact, already first-order dependence logic extended with full intuitionistic negation alone is able  
%including intuitionistic negation alone in first-order dependence logic enables the logic 
 to express finiteness of  the model (via the intuitionistic negation $\neg\phi_\infty$ of the dependence logic sentence $\phi_\infty$ that expresses infinity\footnote{To be more precise, the formula $\phi_\infty$ is defined as $\exists z\forall x\exists y(\dep(x,y)\wedge \dep(y,x)\wedge y\neq z)$, which is true in a first-order model $M$ iff there is an injective function $f:M\to M$ that is not onto, meaning that $M$ is infinite. See e.g., \cite[Section 4.2]{Van07dl} for details.}), and thus compactness fails for the logic.

However, the situation changes significantly when we restrict our attention to propositional and modal team-based logics. Many noteworthy propositional and modal team-based logics are shown to be expressively complete with respect to certain closure properties (such as downward and union closure)
 \cite{HLSV14,HellaStumpf15,Yang2022,VY_PD,YangVaananen:17PT}. The intuitionistic negation preserves most of these closure properties (as $\neg\phi$ is always flat). Consequently,  extending these logics with full intuionistic negation does not increase their expressive power. It is thus both safe 
and reasonable %natural 
 to consider full negation in propositional and modal settings.

%the non-downward closed cases or for $\vvee$ and thus for many of these logics anyway.

 %$\neg\neg\dep(\vec{x},y)$ is no longer equivalent to $\dep(\vec{x},y)$, 

%The only thing that is lost is the connection to game-theoretic semantic, as we will now have 

%$\neg\neg\dep(\vec{p},q)\equiv \dep(\vec{p},q)$ according to game-theoretic semantics, whereas $\neg\neg\dep(\vec{p},q)\equiv \neg\bot\equiv \top$ in our case. But there is no known reasonable game-theoretic semantics for the non-downward closed cases or for $\vvee$ and thus for many of these logics anyway.

%In this sense it is both safe and natural to consider full negation in these contexts. The intuitionistic negation is flat, and thus would preserve all the known closure properties. (On the contrary classical negation will destroy these properties.) This negation is also conservative as $\neg\alpha\equiv\dneg\alpha$ for classical $\alpha$. This is what we will do in this paper. We explore logics with full negation. All preliminaries will be provided in Section 2.

Moreover, incorporating this full negation does not complicate the axiomatization of propositional team-based logics with downward closure property, as we will demonstrate in this paper.
%such as classical propositional logic extended with inquisitive disjunction ($\PL(\vvee)$). 
We %demonstrate 
show in Section 4 that for classical propositional logic extended with global disjunction $\vvee$ ($\PL(\vvee)$) with full negation, a complete natural deduction system  can be obtained by adapting the existing system of $\PL(\vvee)$ with restricted negation  in a natural way.  In Section 4 we also revisit the known systems for $\PL(\vvee)$ \cite{Ciardelli2015,CiardelliIemhoffYang2020,VY_PD} in the literature and present a simpler, folklore-based alternative. In addition, we  document several  folklore observations regarding the system, including the fact the Split axiom of inquisitive logic is derivable from a distributive law in the system of $\PL(\vvee)$.%\todo{provides an axiomatization of the downward closed logic $\PL(\vvee)$ based on the known ones for $\PL(\vvee)$ with restricted negation.}

%$\neg\phi$ is defined (in the usual way) as $\phi\to\bot$, and $\dneg\alpha$ is indeed equivalent to $\alpha\to\bot$ for classical $\alpha$. (The intuitionistic implication was introduced in the context of dependence logic in \cite{AbVan09}.) However, dependence logics were first introduced without implication, and they later often don't include such an implication in the syntax for a reason. Namely that first-order dependence logic with intuitionistic implication is expressively too strong, as strong as full second-order logic \cite{Yang2011}, and thus many nice properties are lost. In fact, already adding intuitionistic negation to dependence logic will make the logic able to express finiteness (using the formula $\neg\phi_\infty$) and thus compactness is lost.

Closely related to the lack of appropriate full negation in the typical formulation of team-based logics is the challenge of defining the complement of a given property. %an expressible property within the logic. 
This is particularly relevant in the proof of the expressive completeness of propositional dependence logic ($\PL(\dep(\cdot))$) as given in \cite{VY_PD}. Similar expressive completeness result for $\PL(\vvee)$ is %or other propositional team logics with the global disjunction $\vvee$ are 
proved in  \cite{VY_PD} by expressing a relevant team property using a formula $\bigvvee_{t\in T}\chi_t$ in disjunctive normal form, where each disjunct $\chi_t$ defines a specific (positive) property that that is easily expressible in classical logic. %, which is easily expressible in the classical language. 
In contrast, for  $\PL(\dep(\cdot))$ without the global disjunction, a conjunctive formula $\bigwedge_{t\in T}\rho_t$ is used as the defining formula in the proof, where each conjunct $\rho_t$ defines the complement of a certain property $P$. One might want to define $\overline{P}$ using the negation of a formula $\phi$ that defines $P$, but the appropriate negation in this context is the classical negation $\sim$, which is lacking in the logic's language. Adding such a classical negation $\sim$ would break the downward closure of $\PL(\dep(\cdot))$ and resulting logic has quite different properties (see \cite{Luck18,Puncochar15} for studies on this and related logics). 
Let us also remark that negated formulas  $\neg\phi$ with our newly introduced full intuitionistic negation do not generally contribute to expressing complemented properties, as $\neg\phi$ is always flat and can thus only define flat properties. Nevertheless,  the property $\overline{P}$  can  be expressed in $\PL(\dep(\cdot))$ by other means.  

In Section 3, we review  these results from the literature and highlight how complemented properties are expressed without directly using negation in $\PL(\dep(\cdot))$. This further explains why full negation has often been excluded from the typical formulations of propositional dependence logic and its variants---the versions without full negation are already sufficiently expressive, after all.

It is worth noting that in the first-order setting, as expected, incorporating classical negation $\sim$ into first-order dependence logic dramatically  increases the expressive power of the logic to full second-order logic too \cite{KontinenNurmi2011}.  %The work \cite{Yang_neg18} explores 
Alternative ways to express complemented properties without increasing the expressive power in the first-order setting have been explored in \cite{Yang_neg18}. 

In Section 3, we also revisit the disjunctive and conjunctive normal forms for  propositional team-based logics. %as well as for the usual classical logic. 
We illustrate how the normal form approach for obtaining expressive completeness results for these logics aligns with the standard normal form technique used in the usual (single valuation-based) classical logic. 
%Inspired by this connection, %in a similar manner, 
Building on these insights, we also prove a new result:  classical propositional logic extended with both dependence and inclusion atoms ($\PL(\dep(\cdot),\subseteq)$) is expressively complete for the class of all team properties that contain the empty team. This logic  has so far been largely overlooked in the literature. Our result corresponds well with the known result for the first-order variant of $\PL(\dep(\cdot),\subseteq)$, which is also expressively complete, in the sense that it captures all properties definable in existential second-order logic \cite{Pietro_I/E}. Our proof also offers yet another way to express complemented properties without directly using negation in propositional team-based logics.

Overall,  in this paper we aim to demonstrate that adding full (intuitionistic) negation to propositional team-based logics does not introduce significant technical complications, especially for downward closed logics. The obstacle to incorporating full negation is, in a sense, more conceptual than technical. For instance, one remaining challenge  is how to meaningfully interpret negated dependence or inclusion atoms. Nevertheless, %To this, our response is, 
%at least 
%We emphasize though that 
there are interesting well-behaved propositional team-based  logics that do not contain these dependency atoms in the syntax, such as $\PL(\vvee)$, which can accommodate full negation well. 
We leave the exploration of these and other issues to future research, which is discussed in Section 5.

\section{Preliminaries}

We first fix our notations for the language of classical propositional logic, denoted as \PL. Let \Prop be a (countable) set of propositional variables. We define formulas of \PL, also called {\em classical formulas}, recursively as
\[\alpha::= p\mid\bot\mid \neg \alpha\mid \alpha\wedge\alpha\mid\alpha\vee\alpha\]
where $p\in \Prop$. The standard classical logic evaluates formulas over {\em valuations}, which are functions $v:\Prop\to\{0,1\}$. The satisfaction relation $v\models\alpha$ is defined inductively as usual, namely, e.g., $v\models p$ iff $v(p)=1$,  %$v\models\bot$ never holds, 
$v\models\neg\alpha$ iff $v\not\models\alpha$, etc.

In this paper we study propositional logics with {\em team semantics}. One such logic, denoted as $\PL(\vvee)$, has the language \PL extended with {\em global disjunction}  $\vvee$\footnote{The disjunction $\vvee$ is also known in the literature as {\em intuitionistic disjunction}, {\em inquisitive disjunction}, and sometimes also {\em Boolean disjunction} or {\em classical disjunction}. In this paper we choose to adopt a more neutral term: {\em global disjunction}.} (in contrast, the disjunction of \PL is referred to as {\em local disjunction}). To be more precise, formulas of $\PL(\vvee)$ are defined recursively as
\[\phi::=p\mid \bot\mid\neg\phi\mid \phi\wedge\psi\mid\phi\vee\psi\mid\phi\vvee\psi\]
where $p\in \Prop$. Unlike most literature on team semantics where negation is restricted to atomic or classical formulas only, in this paper we allow negation to occur in front of arbitrary formulas. %As we shall see, this more relaxed syntax 

Formulas of $\PL(\vvee)$ are evaluated over (propositional) teams. A {\em team} is a set $t$ of valuations $v:\Prop\to\{0,1\}$. We define the satisfaction relation $t\models\phi$ inductively as:
\begin{itemize}
\item $t\models p$ ~~iff~~ $v\models p$ for all $v\in t$;
\item $t\models\bot$ ~~iff~~ $t=\emptyset$
\item $t\models\neg\phi$ ~~iff~~ $\{v\}\not\models\phi$ for all $v\in t$;
\item $t\models\phi\wedge\psi$ ~~iff~~ $t\models\phi$ and $t\models\psi$
\item $t\models\phi\vee\psi$ ~~iff~~ there exist $r,s\subseteq t$ such that $t=r\cup s$, $r\models\phi$ and $s\models\psi$;
\item $t\models\phi\vvee\psi$ ~~iff~~ $t\models\phi$ or $t\models\psi$
\end{itemize}
We write $\phi\models\psi$ if $t\models\phi$ implies $t\models\psi$ holds for all teams $t$. If both $\phi\models\psi$ and $\psi\models\phi$, we write $\phi\equiv\psi$ and say that $\phi$ and $\psi$ are (semantically) equivalent.

The above semantic clause for the full negation $\neg\phi$ is clearly compositional. 
%, and it 
%coincides with that of intuitionistic negation $\phi\to\bot$, where $\to$ is the intuitionistic implication defined as\todo{cite something}
%\begin{itemize}
%\item $t\models\phi\vvee\psi$ ~~iff~~ for any $s\subseteq t$, $s\models\phi$ implies $s\models\psi$.
%\end{itemize}
 We will make  further remarks about the full negation at the end of this section. Let us now focus on recalling basic properties of the logic.

It is easy to verify by induction that arbitrary formulas $\phi$ in $\PL(\vvee)$ satisfy the empty team and downward closure property, defined as:
\begin{description}
\item[Empty team property:] $\emptyset\models\phi$ always holds.
\item[Downward closure:] If $t\models\phi$, then $s\models\phi$ for all $s\subseteq t$.
\end{description}
Formulas in the classical language $\alpha$ satisfy, additionally, union closure and flatness property:
\begin{description}
\item[Union closure:] If $t_i\models\alpha$ for all $i\in I\neq\emptyset$, then $\bigcup_{i\in I}t_i\models\alpha$.
\item[Flatness:] $t\models\alpha$ iff $\{v\}\models\alpha$ for all $v\in t$.
\end{description}
One can easily show that a formula $\phi$ is flat iff it satisfies downward closure, union closure and the empty team property\footnote{Even though this fact is immediate, in some earlier publications on team semantics (such as the author's joint paper \cite{YangVaananen:17PT}), the empty team property was mistakenly forgotten. Neither did they include the restriction that in the definition of union closure, the index set $I$ should be assumed nonempty, for otherwise union closure would (unintentionally) imply the empty team property. The author would like to thank Fausto Barbero for bringing this issue to her attention several years ago.}.

An easy inductive proof shows that for any classical formula $\alpha$,
\begin{equation}\label{flat_singleton_singleval}
\{v\}\models\alpha\iff v\models\alpha.
\end{equation}
We will thus use the above two types of satisfaction relations for classical formulas interchangeably.
%sloppy about the notation for satisfaction relation we use for classical formulas.

%%In fact, flatness is the characteristic property for classical formulas, in the sense of the following proposition (see, e.g., \cite{YangVaananen:17PT}, or \Cref{pl-exp-comp} in the sequel):
%%\begin{prop}
%%A formula is flat iff it is equivalent to a classical formula.\todo{keep or not?}
%%\end{prop}

%Throughout the paper, we reserve the first three Greek letters $\alpha,\beta,\gamma$  (with or without subscripts)  for classical formulas.  

Non-classical formulas are in general not flat. A typical example of a non-flat non-classical formula is $p\vvee\neg p$, which is trivially true on all singleton teams (i.e., teams of the form $\{v\}$) but not true, for example, on the team $\{v_0,v_1\}$ with $v_0(p)=0$ and $v_1(p)=1$. Yet, the negation $\neg(p\vvee\neg p)$ (being still non-classical) is flat. In fact, it follows immediately from the definition that all negated formulas $\neg\phi$ are flat.

%More generally, all {\em Harrop formulas} are  flat. The class of {\em Harrop formulas} is defined syntactically  as
%\[\delta::=p\mid \bot\mid \neg\phi\mid\delta\wedge\delta\mid \delta\vee\delta\]
%where $\phi$ ranges over arbitrary $\PL(\vvee)$ formula. For instance, all classical formulas, $\neg (p\vvee \neg p)$, and $\neg (\phi\vvee\psi)\vee (q\wedge r)$ are  Harrop formulas, whereas $p\vvee \neg p$ is not a Harrop formula. Throughout the paper, we use $\delta,\eta$ (with or without subscripts) to stand for Harrop formulas.  %Since negated formulas are  flat, Harrop formulas are flat.
%Harrop formulas are studied originally in the literature on intuitionistic logic \cite{Harrop56}. They are considered by Pun\u{c}och\'{a}\u{r} \cite{FergusonPuncochar2013,Puncochar2016} also in the context of inquisitive logic and team-based logics. Both contexts consider Harrop formulas in languages with intuitionistic implication, which is not available in $\PL(\vvee)$.\todo{defined implication?} Harrop formulas are of particular interest in our context as well, as they are all easily verified (by induction) to be flat, even though they can contain the global disjunction $\vvee$ (and are thus not classical).\todo{for Harrop, $\neg\neg\delta\equiv\delta$}\todo{Move Harrop discussion to the proof theory section}
%may not be all classical (as they may contain the global disjunction $\vvee$).\todo{one can also replace $\vvee$ with $\vee$ in a Harrop formula}

%Even though Harrop formulas may not be classical, it is easy to verify by induction that they are all flat.

We also consider dependence and inclusion atoms, which are strings of the form $\dep(p_1\dots p_n,q)$ and $a_1\dots a_n\subseteq b_1\dots b_n$ (where $a_i,b_i\in \Prop\cup\{\bot,\top\}$), respectively. Their semantics are defined as
\begin{itemize}
\item $t\models\dep(p_1\dots p_n,q)$ ~~iff~~ $v(p_i)=u(p_i)$ for all $1\leq i\leq n$ implies $v(q)=u(q)$.
\item $t\models a_1\dots a_n\subseteq b_1\dots b_n$ ~~iff~~ for all $v\in t$, there exists $u\in t$ such that  for all $1\leq i\leq n$, $v(a_i)=u(b_i)$.
\end{itemize}
Two simple types of these atoms will play a central role in our discussions: Dependence atoms $\dep(\langle\rangle,p)$ with the first component the empty sequence, which we also denote as $\dep(p)$ and refer to as {\em constancy atoms}; inclusion atoms of the form $x_1\dots x_n\subseteq p_1\dots p_n$ with $x_1,\dots,x_n\in\{\top,\bot\}$, which we refer to as {\em primitive inclusion atoms}. The semantic clauses of these simple atoms reduce to
\begin{itemize}
\item $t\models\dep(p)$ ~~iff~~ for all $v,u\in t$, $v(p)=u(p)$ (namely, $p$ has a constant value in the team $t$).
\item $t\models x_1\dots x_n\subseteq p_1\dots p_n$ ~~iff~~ there exists $v\in t$ such that  for all $1\leq i\leq n$, $v(p_i)=\underline{x_i}$, where $\underline{\top}=1$ and $\underline{\bot}=0$ (namely, the sequence  $\underline{x_1}\dots \underline{x_n}$ of values for $p_1\dots p_n$ is present in the team $t$).
\end{itemize}
%Intuitively, the constancy atom $\dep(p)$ states that the propositional variable $p$ has a constant value in the team in questions, and the primitive inclusion atom, say, $\bot\top\subseteq pq$, states that the value $01$ 

%$\dep(\delta_1\dots\delta_n,\eta)$ and $\delta_1\dots,\delta_n\subseteq \eta_1\dots\eta_n$, respectively, where $\delta_i$ and $\eta_i$ are Harrop formulas. Their semantics are defined as
%\begin{itemize}
%\item $t\models\dep(\delta_1\dots\delta_n,\eta)$ ~~iff~~ $v\sim_{\{\delta_1,\dots,\delta_n\}}u$ implies $v\sim_{\{\eta\}}u$, where $v\sim_\Delta u$ is defined as  for all $\delta\in \Delta$, 
%\[\{v\}\models\delta \iff \{u\}\models\delta.\]
%\item $t\models\delta_1\dots\delta_n\subseteq \eta_1\dots\eta_n$ ~~iff~~ for all $v\in t$, there exists $u\in t$ such that  for all $1\leq i\leq n$,
%\[\{v\}\models\delta_i \iff \{u\}\models\eta_i.\]
%\end{itemize}
%Note that we allow more relaxed syntax for dependence and inclusion atoms than the ones for the so-called {\em extended dependence} and {\em extended inclusion} atoms \cite{}\todo{references}, where the arguments $\delta_1,\dots,\delta_n,\eta,\eta_1\dots\eta_n$ are all assumed to be classical formulas. %(which are special cases of Harrop formulas).

%in case all the Harrop formulas involved are classical formulas, by (\ref{flat_singleton_singleval}), the above definitions coincide with the standard definitions in the literature.\todo{rephrase}

%We form extensions of $\PL$  by incorporating these atoms and global disjunction into the language. 
We write $\PL(\dep(\cdot))$ for $\PL$ extended with dependence atoms, known as {\em propositional dependence logic},  $\PL(\subseteq)$ for $\PL$ extended with inclusion atoms, known as {\em propositional inclusion logic}, $\PL(\subseteq,\vvee)$ for $\PL$ extended with  inclusion atoms and global disjunction, and $\PL(\dep(\cdot),\subseteq)$ for $\PL$ extended with both dependence and inclusion atoms. To be more precise,  formulas of, e.g., $\PL(\dep(\cdot))$ are defined recursively as
\[\phi::=p\mid \bot\mid\dep(p_1\dots p_n,q)\mid
 %\dep(\delta_1\dots\delta_n,\eta)\mid
 \neg\phi\mid \phi\wedge\psi\mid\phi\vee\psi\mid\]
where $p\in \Prop$. %and $\delta_1,\dots,\delta_n,\eta$ are Harrop formulas (which, in this case, are all classical formulas). %\todo{rephrase:}
%On the other hand, in the logic $\PL(\subseteq,\vvee)$, Harrop formulas are not necessarily classical.
It is easy to verify that formulas in $\PL(\dep(\cdot))$ are downward closed, formulas in $\PL(\subseteq)$ are union closed, and formulas in all these four new logics have the empty team property. 
%All these four logics as well as $\PL(\vvee)$ are expressively complete with respect to 

%Define the syntax so that negation is allowed, $\neg\phi$, and even $\dep(\alpha)$ where $\alpha$ is classical or negated (called Harrop formula -- define it!) \todo{cite Vit}

We conclude this section with a few additional remarks on the full negation we defined.  %The reader may have already observed that the above 
Our semantic clause for the full negation $\neg\phi$ %is clearly compositional, and it 
%is 
actually corresponds to that of
the familiar intuitionistic negation $\phi\to\bot$, where $\to$ denotes the {\em intuitionistic implication} from the literature (see e.g.,  \cite{AbVan09,InquiLog}) and is defined as%\todo{cite something}
\begin{itemize}
\item $t\models\phi\to\psi$ ~~iff~~ for any $s\subseteq t$, $s\models\phi$ implies $s\models\psi$.
\end{itemize}
%Intuitionistic implication plays an important role in some other team-based logics, such as inquisitive logic \cite{InquiLog}.
While we, in this paper, do not include the intuitionistic implication in our syntax, adding it to, e.g., $\PL(\vvee)$ does not increase the expressive power of the logic, as was proved essentially in \cite{VY_PD} and as we will review in \Cref{exp-comp-plv}. (The intuitionistic implication is, however, not uniformly definable in $\PL(\vvee)$, by the arugments in \cite{Yang15}.\footnote{Without going into too much detail, we note that any context formed by the connectives in $\PL(\vvee)$ is monotone, while the formula $\neg\phi$ is not. See \cite{Yang15} for the definition of a context and a monotone context.})

As expected for intuitionistic negation, the excluded middle law fails with respect to local or global disjunction in the downward closed logics $\PL(\vvee)$ and $\PL(\dep(\cdot))$, namely neither $\phi\vee\neg\phi$ nor $\phi\vvee\neg\phi$ is valid in general for these logics. To see why, for the logic $\PL(\vvee)$, the two formulas $ (p\vvee\neg p)\vee\neg(p\vvee\neg p)$ and  $p\vvee\neg p$ are both invalid. %and thus serve as counter-examples. 
For $\PL(\dep(\cdot))$, our counter-example is $\dep(p)\vee\neg\dep(p)$, where note that $\neg\dep(p)\equiv \bot$.

Over union closed team-based logics, $\phi\vee\neg \phi$ is actually valid. To see why, consider any team $t$, and let $r=\{v\in t\mid \{v\}\models\phi\}$ and $s=\{v\in t\mid \{v\}\not\models\phi\}$. Since $\phi$ is closed under unions, $r\models\phi$ and $s\models\neg\phi$, which give $t\models\phi\vee\neg\phi$.

\section{Expressive completeness via normal form and ways to negate% (revisited)
}

%\todo{discuss locality}

All  five of the team-based logics we defined are expressively complete. The expressive completeness result for the downward closed logics $\PL(\vvee)$ and $\PL(\dep(\cdot))$ with atomic negation only was proved in  \cite{VY_PD}, while for $\PL(\subseteq)$ with restricted negation, it was proved in \cite{Yang2022}, and the result for $\PL(\subseteq,\vvee)$ with restricted negation is folklore. The expressive completeness of these four logics with full negation %(as we defined) 
follows immediately.
In this section, we review the role of disjunctive and conjunctive normal forms in obtaining the expressive completeness results for these  logics. Our goal is to highlight how this approach aligns with  the standard techniques used in the usual (single valuation-based) classical logic. We also discuss how the complement of a property was expressed in %propositional dependence logic %(with no classical negation in the language) 
 $\PL(\dep(\cdot))$ without directly using negation,
 %these logics with no full negation
  and point out that our more relaxed full negation does not readily offer new methods for capturing these complemented properties.
%Our goal is to highlight that this approach actually resembles the similar standard approach taken in the usual (single valuation-based) classical logic. %We also review how classical negation plays a role in this context. In particular, 
%We also highlight how the complement of a property was expressed in %propositional dependence logic %(with no classical negation in the language) 
% these logics with no full negation. We point out that our more relaxed full negation actually does not readily offer new methods for capturing these negated properties.
Finally, building on the insights of this alternative way of expressing complemented properties, we prove a new result that the logic $\PL(\dep(\cdot),\subseteq)$ is expressively complete, by using conjunctive normal form and by expressing complemented properties in yet another way. %We also give the proof of a folklore result that the logic $\PL(\subseteq,\vvee)$ is expressively complete.

%\todo{say that these results are proved essentially in the literature. Our new negation does not affect anything}

%\todo{All of the five logics are expressively complete}

\subsection{Classical logic and normal forms}\label{sec:pl}

It will be instructive for our discussion to first recall the standard expressive completeness result for classical logic, along with how it can be proved by constructing formulas in disjunctive and conjunctive normal forms.

%Let us spend some words on recalling the standard expressive completeness result for classical logic. We will also recall how this can be proved by constructing formulas in disjunctive and conjunctive normal form.

Throughout Session 3, we fix a set $N=\{p_1,\dots,p_n\}$ of $n$ propositional variables. All formulas considered in this section are in $N$ and all teams are $N$-teams. For any $\PL$-formula $\alpha$ with propositional variables from $N$, consider the set
\[\Arrowvert\alpha\Arrowvert=\{v\in 2^N: v\models \alpha\}\]
of valuations that satisfies $\alpha$. We shall call the set  $\Arrowvert\alpha\Arrowvert\subseteq 2^N$  the {\em property characterised by $\alpha$} (over $N$).
%Clearly, $\Arrowvert\alpha\Arrowvert$ is a subset of $2^N$ (and can thus also be viewed as a team). 
We now review the result that classical formula is {\em expressively complete}, in the sense that for any property $t\subseteq 2^N$ over $N$, classical logic has a {\em characteristic formula} $\chi_t$ satisfying $\Arrowvert \chi_t\Arrowvert=t$, or, in other words, for any valuation $v\in 2^N$,
\begin{equation}\label{cl_char_eqv}
v\models \chi_t\iff v\in t.
\end{equation}
Let us examine two ``brute-force'' methods to construct such a formula. These are based on two alternative interpretations of the equivalence (\ref{cl_char_eqv}), leading to a disjunctive normal form and a conjunctive normal form for classical logic.

We first observe that for an arbitrary single valuation $v\in 2^N$, the classical formula
\begin{equation*}\label{cl_v_char}
\chi_v=p_1^{v(p_1)}\wedge\dots\wedge p_n^{v(p_n)},~~\text{ where }p_i^1:=p_i\text{ and }p_i^0:=\neg p_i,
\end{equation*}
defines the valuation, in the sense that for any valuation $u\in 2^N$,
\begin{equation}\label{cl_v_char_prop}
u\models\chi_v\iff u=v.
\end{equation}
%\todo{remark that one can also just take restrictions, instead of fixing $N$.}
For any property $t\subseteq 2^N$ (which is also a team), we now form the characteristic formula 
\begin{equation}\label{cl_char_eqv_disj}
\chi_t=\bigvee_{v\in t}\chi_v
\end{equation}
by taking disjunctions of the $\chi_v$ formulas. Clearly, the $\chi_t$ formula in disjunctive normal form satisfies (\ref{cl_char_eqv}), as it expresses the property that
\begin{center}
``the valuation $v$ in question is equal to some valuation in $t$''
\end{center}

%Now, our characteristic formula $\chi_t$ that would satisfy (\ref{cl_char_eqv}) should then state 
%\begin{center}
%``the valuation $v$ in question is equal to some valuation in $t$."
%\end{center}
%This idea can be captured by using the disjunction connective through defining
%\begin{equation}\label{cl_char_eqv_disj}
%\chi_t=\bigvee_{v\in t}\chi_v.
%\end{equation}
%The reader can see that the formula $\chi_t$ is a disjunction of conjunctions of literals, and thus in {\em disjunctive normal form}.

There is also an indirect, alternative way to interpret the equivalence in (\ref{cl_char_eqv}), namely to let $\chi_t$ express
\begin{center}
``the valuation $v$ in question is not equal to any valuation in the complement of $t$ (and thus  has to be in $t$).''
\end{center}
This idea can be captured by using conjunction and negation  through defining
\begin{equation}\label{cl_char_eqv_conj}
\chi_t=\bigwedge_{v\in \overline{t}}\neg \chi_v
\end{equation}
where note that the set $\overline{t}=2^N\setminus t$ is finite. 
The formula $\chi_t$ can be easily turned into an equivalent formula
\(\bigwedge_{v\in \overline{t}}(p_1^{1-v(p_1)}\vee\dots\vee p_n^{1-v(p_n)})\)
 in {\em conjunctive normal form} 
by applying the De Morgan laws  and double negation elimination law to each formula $\neg\chi_v$. Let us emphasise that (classical) negation plays an important role in the conjunctive characteristic formula $\chi_t$ in (\ref{cl_char_eqv_conj}). %\todo{negation is always flat in team logics}

%\todo{emphasize the role of negation - in this language that has negation - later in our logics, negation is not always available and its role is also different.}

%Given a valuation $v\in 2^N$, classical logic 
%\[\chi_t=\bigvee_{v\in t}(p_1^{v(1)}\wedge\dots\wedge p_n^{v(n)})\] 
%since clearly, for any valuation $t\in 2^N$,
%\[v\models \chi_t\iff v\in t\]

%\todo{check some literature about normal forms in classical logic}

\subsection{Team-based logics: downward closed case}\label{sec:dw-exp}
We now turn to reviewing the proofs given in the literature for  the expressive completeness of the downward closed team-based logics $\PL(\vvee)$ and $\PL(\dep(\cdot))$. %, which actually %We highlight how these proofs 
%reply on  essentially similar ideas to the standard arguments in classical logic.
%In this section we prove that the logic $\PL(\dep(\cdot),\subseteq)$ is expressively complete via the method of normal form. 

Define a {\em team property} (over $N$) as a set $T\subseteq \wp(2^N)$ of teams. The {\em team property characterised by a formula $\phi$} (over $N$) is defined as the set
\[\llbracket\phi \rrbracket=\{t\subseteq 2^N: t\models\phi\}\]
of $N$-teams that satisfies $\phi$. As an example, the team property characterised by a classical formula $\alpha$ can be easily computed from the standard (single-valuation) property $\Arrowvert \alpha\Arrowvert$ characterised by the formula:

%Consider the set
%\[\llbracket\phi \rrbracket=\{t\subseteq 2^N: t\models\phi\}\]
%of $N$-teams that satisfies $\phi$.  We shall call the set  $\llbracket\phi \rrbracket\subseteq \wp(2^N)$  the {\em team property characterised by $\phi$} (over $N$).

\begin{lem}\label{pl_pw_truth}%\todo{combine this with \Cref{theta_frm_prop}?}
%For any classical formula $\alpha$, we have that 
\(\llbracket\alpha \rrbracket=\wp(\Arrowvert \alpha\Arrowvert)\) for  any classical formula $\alpha$.
\end{lem}
\begin{proof}
For any  $t\in \llbracket\alpha \rrbracket$, we have $t\models\alpha$. Thus, $v\models\alpha$ for any $v\in t$, by flatness of $\alpha$. It then follows that $t\subseteq \Arrowvert \alpha\Arrowvert$, i.e., $t\in \wp(\Arrowvert \alpha\Arrowvert)$.

Conversely, for any nonempty $t\in \wp(\Arrowvert \alpha\Arrowvert)$, i.e., $\emptyset\neq t\subseteq \Arrowvert \alpha\Arrowvert$, we have $t\models\alpha$ by flatness of $\alpha$ since $v\models\alpha$ for each $v\in t$. The fact $\emptyset\models\alpha$ follows from the empty team property of $\alpha$.
\qed\end{proof}

%Define 
%\[\Theta_t=\bigvee_{v\in t}(p_1^{v(1)}\wedge\dots\wedge p_n^{v(n)}).\]

The team property $\wp(\Arrowvert \alpha\Arrowvert)$ is clearly {\em flat}, meaning (naturally) that $t\in \wp(\Arrowvert \alpha\Arrowvert)$ iff $\{v\}\in \wp(\Arrowvert \alpha\Arrowvert)$ for all $v\in t$. We can define  a team property $T$ being {\em downward closed} or {\em union closed} in a similar way.
%Similarly, we can speak about a team property $T$ being downward closed or union closed.

Let $\mathcal{C}$ be a class of team properties. We say that a team-based logic $\LL$ is {\em expressively complete for $\mathcal{C}$} if the following two conditions hold:
\begin{itemize}
\item for every $\LL$-formula $\phi$, the team property $\llbracket\phi \rrbracket$  belongs to $\mathcal{C}$, i.e., $\llbracket\phi \rrbracket\in \mathcal{C}$;
\item for every team property $T\in\mathcal{C}$, there exists an $\LL$-formula $\phi$ such that $\llbracket\phi \rrbracket=T$.
\end{itemize}

%\begin{thm}
%\PL is expressively complete for the class of flat team properties.
%\end{thm}
%\begin{proof}
%\end{proof}

\begin{thm}\label{exp-comp-plv}
The logics $\PL(\vvee)$ and $\PL(\dep(\cdot))$ are both expressively complete for the class of nonempty and downward closed team properties. %, where a team property $T$ being downward closed means (naturally) that $s\subseteq t\in T$ implies $s\in T$.
\end{thm}

%\todo{cite Ivano's master thesis and say that there it is shown essentially that this $\PL(\vvee)$ with our full negation is expressively complete}
%The downward closed logics $\PL(\vvee)$ and $\PL(\dep(\cdot))$ were proved to be expressively complete with respect to nonempty downward closed team properties 
The above theorem are proved essentially already in  \cite{VY_PD} (where, however, arbitrary negation was not allowed in the syntax). We now revisit these proofs and demonstrate how the standard approach of constructing formulas in disjunctive and conjunctive normal form from the previous subsection can be extended to the context of team-based logics, with a twist regarding negation.

%We now review these proofs and illustrate how the standard idea of constructing formulas in disjunctive and conjunctive normal form from the previous subsection can be generalised to the team-based logics case, with a twist concerning negation.

The first condition in the expressive completeness definition for the two logics $\PL(\vvee)$ and $\PL(\dep(\cdot))$ follows easily from the fact that formulas $\phi$ in the  logics are downward closed and have the empty team property (thus $\emptyset\in\llbracket\phi \rrbracket\neq\emptyset$).  For the non-trivial second condition, we need to characterise every nonempty downward closed team property $T$ using a formula $\chi_T$ in the two logics; that is to show that 
\begin{equation}\label{team_char_eqv}
s\models \chi_T\iff s\in T.
\end{equation}
This equivalence is the team semantics counterpart to equivalence (\ref{cl_char_eqv}) from the previous subsection, where  the classical formula $\chi_t$ was used as a characteristic formula. %, and we shall prove it using a similar approach with disjunctive and conjunctive normal forms.

 As a starting point of our proof, observe that on the team level, the classical characteristic formula $\chi_t$ (given by either (\ref{cl_char_eqv_disj}) or (\ref{cl_char_eqv_conj})) defines the team $t$ up to its subteams. %, or it  characterises the team property $\wp(t)$. 
 
% the team property $\wp(t)$ for any team $t$ (over $N$) can be characterised by the  classical characteristic formula $\chi_t$ we discussed in the previous subsection:

%Thus, the classical characteristic formula $\chi_t$ we discussed in the previous subsection characterises the team property $\wp(t)$:
\begin{lem}\label{theta_frm_prop}
For any $N$-team $s$, 
\[s\models\chi_t\iff s\subseteq t.\]
In other words, $\llbracket \chi_t\rrbracket=\wp(t)$.
\end{lem}
\begin{proof}
By \Cref{pl_pw_truth} and (\ref{cl_char_eqv}), $\llbracket\chi_t \rrbracket=\wp(\Arrowvert \chi_t\Arrowvert)=\wp(t)$. %, which gives the required equivalence.
\qed\end{proof}

A corollary that the classical language \PL over team semantics is also expressively complete follows easily:%\todo{used?}
\begin{cor}\label{pl-exp-comp}
\PL is expressively complete for the class of flat team properties. 

Consequently, a formula is flat iff it is equivalent to a classical formula.
\end{cor}
\begin{proof}
For the expressive completeness, the first condition follows from the fact that classical formulas are flat. For the second condition, for any flat team property $T$, observe that $T=\wp(\bigcup T)$. Thus, by \Cref{theta_frm_prop}, we have $\llbracket \chi_{\bigcup T}\rrbracket=\wp(\bigcup T)=T$, where $ \chi_{\bigcup T}$ is a classical formula.

For the second claim in the corollary, classical formulas are flat. Conversely, any flat formula $\phi$ defines a flat property $\llbracket\phi\rrbracket$, which can then be defined by a \PL-formula $\alpha$., i.e., $\llbracket\phi\rrbracket=\llbracket\alpha\rrbracket$. Thus, $\phi\equiv\alpha$.
\qed\end{proof}

Now,  in the logic $\PL(\vvee)$ we define the characteristic formula $\chi_T:=\bigvvee_{t\in T}\chi_t$ as a disjunction of the classical characteristic formulas $\chi_t$. The reader can verify that (\ref{team_char_eqv}) holds for the formula, since $T$ is downward closed and nonempty (thus $\emptyset\in T$); see \cite{VY_PD} for proof details. Essentially, the formula $\chi_T$ states 
\begin{center}
``the team $s$ in question is a subteam of some team $t$ in $T$.''
\end{center}

% can define the characteristic formula $\chi_T$ for the equivalence (\ref{team_char_eqv}) by using the formula $\chi_t$ and the global disjunction:
%
%\begin{thm}\label{exp-comp-plv}
%The logic $\PL(\vvee)$ is expressively complete for the class of nonempty and downward closed team properties. %, where a team property $T$ being downward closed means (naturally) that $s\subseteq t\in T$ implies $s\in T$.
%\end{thm}
%
%\begin{proof}
%To verify the second condition, for any nonempty  downward closed team property $T\subseteq \wp(2^N)$, define $\chi_T=\bigvvee_{t\in T}\chi_t$. By \Cref{theta_frm_prop}, $\llbracket\chi_T \rrbracket=\bigcup_{t\in T}\llbracket\chi_t \rrbracket=\bigcup_{t\in T}\wp(t)=T$, since $T\neq \emptyset$ is downward closed. We refer the reader to  \cite{VY_PD} for proof details.
%\qed\end{proof}

%\todo{explain the (different and related) ideas of disjunctive normal form, and conjunctive normal form}

%\todo{An example of the conjunctive normal form argument from the literature:}
The global disjunction $\vvee$, which plays a crucial role in the above argument, is however not available in the logic $\PL(\dep(\cdot))$. Yet in $\PL(\dep(\cdot))$ we can construct $\chi_T$ as a conjunctive formula $\chi_T=\bigwedge_{t\in \overline{T}}\rho_t$ by using a similar idea to the previous subsection. One naive idea would be to take $\rho_t=\neg \chi_t$, as we did in  (\ref{cl_char_eqv_conj}) in classical logic. However, our negation  $\neg$ does not behave classically on the team level. In particular, $\neg\chi_t$ is always flat, and thus $\bigwedge_{t\in \overline{T}}\neg\chi_t$ would also be flat, which already means that it cannot be the defining formula for any downward closed team property $T$ that is not flat.

%we %actually 
%have
%\[s\models \neg\chi_t\iff s\cap t=\emptyset,\]
%from which  the equivalence (\ref{team_char_eqv}) does not follow\footnote{To see why, consider $T=\{\{v\},\emptyset\}$ for some $v\in 2^N$. For the team $s=\{v\}\in T$, $s\models  \bigwedge_{t\in \overline{T}}\neg\chi_t$ actually fails, because, e.g., for the team  $t_0=\{v,u\}$ where $v\neq u$, since $s\cap t_0\neq\emptyset$, we must have $s\not\models\neg\chi_{t_0}$.
%%. For any $s\in 2^N$ such that $s\models \bigwedge_{t\in \overline{T}}\neg\chi_t$, we must have $s\models \neg\chi_{t_0}$ as $t_0\in \overline{T}$. But then, $s\cap t_0=\emptyset$ implying that $v\notin s$
%}. 

%Recall that in the conjunctive formula $\chi_t$ given by (\ref{cl_char_eqv_conj}), (classical) negation played a crucial  role. However, our negation connective $\neg$ does not behave classically on the team level. In particular, 
%\[s\models \neg\chi_t\iff s\cap t=\emptyset\]

The proof in \cite{VY_PD}  instead builds the formula $\rho_t$ as in the lemma below:

\begin{lem}\label{taneli_frm}
For any nonempty $N$-team $t$, there is a formula  $\rho_t$ in $\PL(\dep(\cdot))$ such that for any $N$-team $s$,
\[s\models \rho_t\iff t\nsubseteq s.\]
\end{lem}
\begin{proof}
The argument is due to Huuskonen. %We only sketch the proof as presented in \cite{VY_PD}. 
Let $|t|=k+1$. Define $\rho_t=\mu_k\vee\chi_{\overline{t}}$, where 
\[\mu_k=\bigvee_{i=1}^k(\dep(p_1)\wedge\dots\wedge\dep(p_n))\] 
and $\overline{t}=2^N\setminus t$. To show that the formula $\rho_t$ is as required, one applies \Cref{theta_frm_prop} for the classical characteristic formula $\chi_{\overline{t}}$ and a crucial observation for the formula $\mu_k$ that for any $N$-team $s$, 
\[s\models\mu_k\iff |s|\leq k.\]
We refer the reader to  \cite{VY_PD} for details.
\qed\end{proof}

Then, the conjunctive formula $\chi_T=\bigwedge_{t\in \overline{T}}\rho_t$ would satisfy (\ref{team_char_eqv}) as it states
 \begin{center}
``the team $s$ in question is not a superteam of any team  in the complement of $T$ (and thus has to be in $T$).''
\end{center}
Note that $T$ was assumed to be downward closed and nonempty, which implies $\emptyset\in T$. Thus, every team $t\in \overline{T}$ is nonempty, and each conjunct $\rho_t$ in the formula $\chi_T$ is guaranteed to exist by  \Cref{taneli_frm}.

%$\neg\chi_t$ simply does not express the property we wanted, and what we need is classical negation (which is not in the language, and adding it will take us outside of downward closed environment -- in $\PL(\dot\sim)$, one can use this trick to get conjunctive normal form (check! and check whether Martin L\"{u}ck has done it)

%\todo{construction site}
%
%
%
%\begin{thm}
%Propositional dependence logic $\PL(\dep(\cdot))$ is expressively complete.
%\end{thm}
%\begin{proof}
%For any nonempty and downward closed team property $T\subseteq \wp(2^N)$ over $N$, define $\chi_T=\bigwedge_{t\in \overline{T}}\rho_t$, where $\overline{T}=\wp(2^N)\setminus T$. It follows from \Cref{taneli_frm} that $\llbracket\chi_T \rrbracket=T$. We refer the reader to  \cite{VY_PD} for proof details.
%\qed\end{proof}

\subsection{The team-based logics with inclusion atoms}

In this section, we focus on the remaining three team-based logics $\PL(\subseteq)$, $\PL(\subseteq,\vvee)$ and $\PL(\dep(\cdot),\subseteq)$ with inclusion atoms. %, which do not have the downward closure property. 

%We review how a similar disjunctive normal form was used in the literature to prove the expressive completeness of $\PL(\subseteq)$, and how this approach can be generalised to obtain the folklore result that $\PL(\subseteq,\vvee)$ is also expressively complete. 
%%and in folklore result to $\PL(\subseteq)$ and $\PL(\subseteq,\vvee)$ to obtain expressive completeness results. 
%Finally, we present a new result that the logic $\PL(\dep(\cdot),\subseteq)$ is expressively complete, using an argument based on a  conjunctive normal form obtained in a similar and yet technically different manner.

 %for the class of all team properties that contain the empty team.

Recall that for any valuation $v\in 2^N$, in the usual classical logic \PL, we have a formula %$\chi_v=p_1^{v(1)}\wedge\dots\wedge p_n^{v(n)}$ 
$\chi_v$ that defines $v$ in the sense of (\ref{cl_v_char_prop}), which, on the team level, can  be rendered as
\[\{u\}\models\chi_v\iff u=v.\]
This same equivalence can also be captured by the following inclusion atom:
\[\iota_v=\underline{v(p_1)}\dots\underline{v(p_n)}\subseteq p_1\dots p_n,\]
where $\underline{v(p_i)}=p_i$ if $v(p_i)=1$, and $\underline{v(p_i)}=\neg p_i$ if $v(p_i)=0$. That is,
\begin{equation*}\label{inc_iota_char_prop}
\{u\}\models \iota_v\iff u=v.
\end{equation*}
Moreover, %with respect to nonempty teams of arbitrary  cardinality, we have that 
for any nonempty team $s\subseteq 2^N$, 
\begin{equation}\label{inc_char_prop}
s\models\iota_v \iff v\in s.
\end{equation}

Now, as in the case for \PL, the formula $\iota_v$ can  be used as a building block to construct characteristic formulas for teams, in conjunctive and disjunctive normal forms. First, consider the formula $\iota_t=\bigwedge_{v\in t}\iota_v$ (due to  \cite{HellaStumpf15,Yang2022}). %for any nonempty team $t$. 
By (\ref{inc_char_prop}), the formula defines the team $t$ up to its super-teams modulo the empty team, in the following sense:
\begin{equation}\label{sup_team_char}
s\models \iota_t\iff t\subseteq s\text{ or }s=\emptyset.
\end{equation}

%We now use $\iota_v$'s as building blocks to form interesting formulas in conjunctive and disjunctive normal form.
%
%The result below is due to \cite{HellaStumpf15,Yang2022}. The following resembles a conjunctive normal form of \PL.
%
%\begin{lem}\label{inc_charac}
%For any nonempty $N$-team $t\neq \emptyset$, there is a formula $\iota_t$ in $\PL(\subseteq)$ such that for any $N$-team $s$,
%\[s\models \iota_t\iff t\subseteq s\text{ or }s=\emptyset.\]
%\end{lem}
%\begin{proof}
%%\todo{A result  in my paper 2022 PU paper, mention also Hella et al.}
%Define 
%\(\iota_t=\bigwedge_{v\in t}\iota_v.\)
%%where $\underline{v(p_i)}=p_i$ if $v(p_i)=1$, and $\underline{v(p_i)}=\neg p_i$ if $v(p_i)=0$. 
%The result then follows from (\ref{inc_char_prop}). 
%%the crucial observation that  for any nonempty team $r\subseteq 2^N$, 
%%\begin{equation}
%%$r\models\underline{v(p_1)}\dots\underline{v(p_n)}\subseteq p_1\dots p_n$ iff $v\in r$.
%%\end{equation}
%We refer the reader to \cite{Yang2022} for proof details.
%\qed\end{proof}

%\todo{remark that the inclusion atom used above is also a description of a valuation, but it does not define it. However, $\{u\}\models \underline{v(p_1)}\dots\underline{v(p_n)}\subseteq p_1\dots p_n$ iff $u=v$}

Recall that the \PL-formula $\chi_t$ defines the team $t$ up to its subteams in the sense of \Cref{theta_frm_prop}. Thus, %for any nonempty team $t$, 
the formula $\chi_t\wedge \iota_t$ defines the team $t$ precisely modulo the empty team, in the following sense:
\begin{equation}\label{fixing_team_mod}
s\models \chi_t\wedge \iota_t\iff s=t\text{ or }s=\emptyset.
\end{equation}

By using this fact, it was shown in  \cite{Yang2022} that propositional inclusion logic $\PL(\subseteq)$ is expressively complete for the class of union closed team properties that contain the empty team. The idea for the nontrivial condition is that for any such %union closed 
team property $T\subseteq \wp(2^N)$, %with $\emptyset\in T$, 
one proves that the formula $\bigvee_{t\in T} (\chi_t\wedge \iota_t)$ in  disjunctive normal form defines $T$, namely, 
%constructs a formula $\bigvee_{t\in T} (\chi_t\wedge \iota_t)$ in  disjunctive normal form and shows 
$T=\llbracket \bigvee_{t\in T} (\chi_t\wedge \iota_t)\rrbracket$, by using the fact that $T$ is union closed and $\emptyset \in T$. See  \cite{Yang2022} for the proof details.

By adopting a different disjunction---global disjunction---in the characteristic formula, we now provide the proof of a %new (or 
folklore %) 
 result that the logic $\PL(\subseteq,\vvee)$ is  expressively complete for a larger class of team properties.
\begin{thm}
The logic $\PL(\subseteq,\vvee)$ is expressively complete for the class of team properties that contain the empty team.
\end{thm}
\begin{proof}
%\todo{A result essentially already obtained in my paper with Matilda, or not? Folklore?}
For the nontrivial condition, for any team property $T\subseteq \wp(2^N)$ with $\emptyset\in T$, define $\chi_T=\bigvvee_{t\in T}(\chi_t\wedge \iota_t)$. 
%By (\ref{fixing_team_mod}), we have $\llbracket \chi_t\wedge \sigma_t\rrbracket=\{\emptyset,t\}$ for each $t\in T$. Thus, 
Now, by (\ref{fixing_team_mod}), we have that for any $N$-team $s$,
\[s\models \chi_T\iff s\models\chi_t\wedge\iota_t\text{ for some }t\in T \iff s\in T \text{ or }s=\emptyset\iff s\in T,\]
namely $\llbracket \chi_T\rrbracket=T$.
\qed\end{proof}

%\todo{see also my PU paper for how to define arbitrary inclusion atoms with the primitive ones}The following uses a disjunctive normal form.

Finally, let us present our new result that the logic $\PL(\dep(\cdot),\subseteq)$ is also expressively complete for the class of team properties that contain the empty team. The global disjunction $\vvee$ that played a crucial role in the proof above is not available in the language $\PL(\dep(\cdot),\subseteq)$. We thus have to take a different approach to construct the characteristic formula. The idea is that we shall define each team property $T$ with $\emptyset\in T$ by using a $\PL(\dep(\cdot),\subseteq)$-formula $\bigwedge_{t\in \overline{T}}\psi_t$ in conjunctive normal form with each $\psi_t$ satisfying %for any $N$-team $s$,
\begin{equation}\label{neg_form_psi}
s\models \psi_t\iff s\neq t.
\end{equation}
Then, as with the classical logic case we reviewed in \Cref{sec:pl}, the conjunctive formula would then state that 
\begin{center}
``the team in question is not equal to any team $t$ in the complement of $T$ (and thus has to be in $T$).'' 
\end{center}
%\begin{equation}\label{philo_char_prop}
%\begin{split}
%\text{
%``the team in question is not equal to any team $t$ in the complement of $T$}\\ 
%\text{(and thus has to be in $T$).'' 
%}
%\end{split}
%\end{equation}

Let us emphasise (again) that defining $\psi_t$ to be the negated formula $\neg(\chi_t\wedge\iota_t)$ does not give the equivalence (\ref{neg_form_psi}) (or equivalently $\llbracket\neg(\chi_t\wedge\iota_t)\rrbracket\neq \wp(2^N)\setminus\{t\}$), since, for instance, the team property $\llbracket\neg(\chi_t\wedge\iota_t)\rrbracket$ is always flat, while $\wp(2^N)\setminus\{t\}$ need not be.

Now, in order to construct the right formula $\psi_t$, we first recall that the formula $\iota_t$ satisfying (\ref{sup_team_char})  was built by taking conjunction of all $\iota_v$'s with $v\in t$. Interestingly, the dual property of (\ref{sup_team_char}) can be defined by a dual formula $\sigma_t=\bigvee_{v\in \bar{t}}\iota_v$, formed by taking (local) disjunction of $\iota_v$'s, as we show in the next lemma: %For any $N$-team $t$, define $\sigma_t=\bigvee_{v\in \bar{t}}\iota_v$.
%\todo{rephrase}

%the (local) disjunction of all these $\iota_v$'s (with $t$ nonempty), $\sigma_t=\bigvee_{v\in \bar{t}}\iota_v$ defines a dual property, as we show in the next lemma:

\begin{lem}\label{sigma_inc}
Given any nonempty $N$-team $t$, %Define $\sigma_t=\bigvee_{v\in \bar{t}}\iota_v$. 
%For any nonempty $N$-team $t\neq \emptyset$, define $\sigma_t=\bigvee_{v\in \bar{t}}\iota_v$. 
we have that  
%there is a formula $\sigma_t$ in $\PL(\subseteq)$ such that 
for any $N$-team $s$,
\[s\models \sigma_t\iff s\nsubseteq t\text{ or }s=\emptyset.\]
\end{lem}
\begin{proof}
%First, for an arbitrary valuation $v\in 2^N$, define 
%\[\iota_v=\underline{v(p_1)}\dots\underline{v(p_n)}\subseteq p_1\dots p_n,\]
%where $\underline{v(p_i)}=p_i$ if $v(p_i)=1$, and $\underline{v(p_i)}=\neg p_i$ if $v(p_i)=0$.
%Observe that for any team $s\subseteq 2^N$, 
%\[s\models\iota_v\iff v\in s.\]
%
%Define 
%\[\sigma_t=\bigvee_{v\in \overline{t}}\underline{v(p_1)}\dots\underline{v(p_n)}\subseteq p_1\dots p_n.\]
%where $\underline{v(p_i)}=p_i$ if $v(p_i)=1$, and $\underline{v(p_i)}=\neg p_i$ if $v(p_i)=0$.
If $s=\emptyset$, then $s\models\sigma_t$ follows from the empty team property of  $\sigma_t$. 
Now assume that $s\neq \emptyset$. %First observe that  for any nonempty team $r\subseteq 2^N$, 
%\[r\models\underline{v(p_1)}\dots\underline{v(p_n)}\subseteq p_1\dots p_n\iff v\in r.\]
Suppose $\emptyset\neq s\nsubseteq t$. Then there exists $v\in s$ such that $v\in \bar{t}$. By (\ref{inc_char_prop}), we obtain $s\models \iota_v$, and thus $s\models \sigma_t $.

Conversely, if $\emptyset\neq s\models \sigma_t$, then $s=\bigcup_{v\in \bar{t}}s_v$ for some $s_v\subseteq s$ and each $s_v\models \iota_v$. At least one such $s_v$ with $v\in \bar{t}$ must be nonempty and satisfies $v\in s_v$ by (\ref{inc_char_prop}). It then follows that $v\in s_v\setminus t$, implying $s\nsubseteq t$.
\qed\end{proof}

Now, we can build the formula $\psi_t$ that satisfies (\ref{neg_form_psi}).

\begin{lem}\label{say_no2team}
Given any nonempty $N$-team $t$, we have that for any $N$-team $s$,
\[s\models \rho_t\vee\sigma_t\iff s\neq t.\]
%In other words, $\llbracket \rho_t\vee\sigma_t\rrbracket=\wp(2^N)\setminus \{t\}$
\end{lem}
\begin{proof}
Suppose $s\neq t$. Then either $t\nsubseteq s$ or $s\nsubseteq t$. In the former case, $s\models \rho_t$ by \Cref{taneli_frm}; in the latter case, $s\models\sigma_t$ by \Cref{sigma_inc}. Thus, in both cases we obtain $s\models\rho_t\vee\sigma_t$ by the empty team property.

Conversely, suppose $s\models \rho_t\vee\sigma_t$. Then $s=s_1\cup s_2$, $s_1\models\rho_t$ and $s_2\models\sigma_t$ for some $s_1,s_2\subseteq s$. If $s_2\neq \emptyset$, then by \Cref{sigma_inc} we have that $s_2\nsubseteq t$, which implies $s\nsubseteq t$, meaning $s\neq t$. If $s_2=\emptyset$, then $s_1=s$. By \Cref{taneli_frm} we have that $t\nsubseteq s_1= s$, meaning $t\neq s$.
\qed\end{proof}

Finally, we derive the expressive completeness of $\PL(\dep(\cdot),\subseteq)$ as an immediate corollary:
\begin{thm}
The logic $\PL(\dep(\cdot),\subseteq)$ is expressively complete for the class of  team properties that contain the empty team.
\end{thm}
\begin{proof}
For the nontrivial condition, for any team property $T\subseteq \wp(2^N)$ with $\emptyset\in T$, define $\chi_T=\bigwedge_{t\in \overline{T}}(\rho_t\vee\sigma_t)$, where note that $\emptyset\notin \overline{T}$. Thus, $\llbracket \chi_T\rrbracket=T$ follows from \Cref{say_no2team}.
%
%By \Cref{say_no2team}, we have that for any $N$-team $s$,
%\[s\models \chi_T\iff s\models \rho_t\vee\sigma_t\text{ for all }t\in \overline{T}\iff s\neq t\text{ for all }t\in \overline{T}\iff s\in T,\]
%namely, $\llbracket \chi_T\rrbracket=T$.
\qed\end{proof}

Clearly the characteristic formulas used in this proof and the other proofs in this section only use atomic negations. Therefore the results in this section also apply to those logics with atomic negation only.

\section{Natural deduction system for $\PL(\vvee)$ and its extensions revisited}

%In this section, we 

%

In this section, we discuss the axiomatization  for the logic $\PL(\vvee)$ and highlight that the natural deduction system for $\PL(\vvee)$ with restricted negation, as known in the literature, remains complete for our version of  $\PL(\vvee)$ with full negation.
 In this sense, the inclusion of full intuitionistic negation in the syntax of $\PL(\vvee)$ does not complicate its axiomatizatin. We also make some remarks about the known proof system   of the extension of $\PL(\vvee)$ with intuitionistic implication.

The first complete natural deduction system for $\PL(\vvee)$ (with only atomic negation) was introduced in \cite{VY_PD}. Later, \cite{Ciardelli2015} introduced a complete natural deduction system for $\PL(\vvee)$ extended with intuitionistic implication $\to$ (with negation defined via the implication as $\neg\phi:=\phi\to\bot$), which combines and simplifies rules in the system of \cite{VY_PD} and the axioms in the Hilbert system of inquisitive logic  \cite{InquiLog}. Such a system was further simplified and generalised to the intuitionistic case in \cite{CiardelliIemhoffYang2020}. It has been folklore that a complete and simpler system for $\PL(\vvee)$ (with negation for classical formulas only) can be obtained by omitting the rules for implication in the system of \cite{CiardelliIemhoffYang2020}. 
We now present a natural deduction system for $\PL(\vvee)$ based on the system in \cite{CiardelliIemhoffYang2020} with  rules for negation adapted  for the (full) intuitionistic negation.

%on a system for $\PL(\vvee)$ extended with intuitionistic implication given in \cite{Ciardelli2015}.

\begin{defn}\label{plvvee_deduction_system}
The system of natural deduction for $\PL(\!\vvee\!)$ consists of the following rules, where the formula $\alpha$ ranges over classical formulas only: %where the formula $\delta$ ranges over Harrop formulas only.
%\todo{source: cite also Ivano, state it as folklore}
\begin{description}
\item[Rule for $\bot$:]\ 

\vspace{2pt}

\begin{tabular}{|C{0.9\linewidth}|}
\hline
\AxiomC{}\noLine\UnaryInfC{$\bot$} \RightLabel{$\bot\textsf{E}$}\UnaryInfC{$\phi$}\noLine\UnaryInfC{}\DisplayProof
\\
\hline
\end{tabular}

\vspace{2pt}

\item[Rules for $\wedge$:]\ 

\vspace{2pt}

%
%\begin{table}[htp]
%\caption{Rules for $\bot,\top,\wedge$}
\begin{tabular}{|C{0.45\linewidth}C{0.45\linewidth}|}
\hline
%\AxiomC{}\noLine\UnaryInfC{$\bot$} \RightLabel{$\bot\textsf{E}$}\UnaryInfC{$\phi$}\noLine\UnaryInfC{}\DisplayProof&\AxiomC{}\noLine\UnaryInfC{} \RightLabel{$\top\textsf{I}$}\UnaryInfC{$\top$}\noLine\UnaryInfC{}\DisplayProof
%\\
\AxiomC{}\noLine\UnaryInfC{$\phi$}\AxiomC{}\noLine\UnaryInfC{$\psi$} \RightLabel{$\wedge\textsf{I}$}\BinaryInfC{$\phi\wedge\psi$}\noLine\UnaryInfC{}\DisplayProof&\AxiomC{}\noLine\UnaryInfC{$\phi\wedge\psi$} \RightLabel{$\wedge\textsf{E}$}\UnaryInfC{$\phi$}\noLine\UnaryInfC{}\DisplayProof~~\AxiomC{}\noLine\UnaryInfC{$\phi\wedge\psi$} \RightLabel{$\wedge\textsf{E}$}\UnaryInfC{$\psi$}\noLine\UnaryInfC{}\DisplayProof\\\hline
\end{tabular}
%\label{pd_rules_bot_top_wedge}
%\end{table}%

\vspace{2pt}

\item[Rules for $\neg$:]\

\vspace{2pt}

%\begin{table}[htp]
%\caption{Rules for $\neg$}
\begin{tabular}{|C{0.94\linewidth}|}
\hline
\AxiomC{}\noLine\UnaryInfC{$[\phi]$}\noLine\UnaryInfC{$\vdots$}\noLine\UnaryInfC{$\bot$} \RightLabel{$\neg\textsf{I}$}\UnaryInfC{$\neg\phi$} \noLine\UnaryInfC{}\DisplayProof
%\quad\quad \AxiomC{$\bot$}\RightLabel{$\mathsf{ex~ falso}$\todo{derivable}}\UnaryInfC{$\phi$} \DisplayProof
\quad\quad\AxiomC{$\phi$}\AxiomC{$\neg\phi$}\RightLabel{$\neg\textsf{E}$}\BinaryInfC{$\psi$} \DisplayProof
 %\quad\quad\AxiomC{$\neg\neg\alpha$}\RightLabel{\nnege\todo{not needed}}\UnaryInfC{$\alpha$} \DisplayProof
\quad\quad \AxiomC{}\noLine\UnaryInfC{$[\neg\alpha]$}\noLine\UnaryInfC{$\vdots$}\noLine\UnaryInfC{$\bot$}\RightLabel{\textsf{RAA}}\UnaryInfC{$\alpha$} \noLine\UnaryInfC{}\DisplayProof\\\hline
\end{tabular}
%\label{pd_rules_neg}
%\end{table}%

\vspace{2pt}

\item[Rules for $\vvee$:]\

\vspace{2pt}

%\begin{table}[htp]
%\caption{Rules for $\vvee$}
\begin{tabular}{|C{0.45\linewidth}C{0.45\linewidth}|}\hline
 \multirow{2}{*}{\AxiomC{$\phi$} \RightLabel{$\vvee\!\textsf{I}$}\UnaryInfC{$\phi\vvee\psi$}\DisplayProof~\AxiomC{$\phi$} \RightLabel{$\vvee\!\textsf{I}$}\UnaryInfC{$\psi\vvee\phi$}\DisplayProof}&\AxiomC{$\phi\vvee\psi$}\AxiomC{}\noLine\UnaryInfC{$[\phi]$}\noLine\UnaryInfC{$\vdots$}\noLine\UnaryInfC{$\chi$} \AxiomC{}\noLine\UnaryInfC{$[\psi]$}\noLine\UnaryInfC{$\vdots$}\noLine\UnaryInfC{$\chi$}\RightLabel{$\vvee\!\textsf{E}$} \TrinaryInfC{$\chi$}\noLine\UnaryInfC{}\DisplayProof
\\\hline
\end{tabular}
%\label{pd_rules_vvee}
%\end{table}%

\vspace{2pt}

\item[Rules for $\vee$:]\

\vspace{2pt}

%\begin{table}[htp]
%\caption{Rules for $\vee$}
\begin{tabular}{|C{0.45\linewidth}C{0.45\linewidth}|}
\hline
 \multirow{2}{*}{\AxiomC{$\phi$} \RightLabel{$\vee\textsf{I}$}\UnaryInfC{$\phi\vee\psi$}\DisplayProof}&\AxiomC{$\phi\vee\psi$}\AxiomC{$[\phi]$}\noLine\UnaryInfC{$\vdots$}\noLine\UnaryInfC{$\alpha$} \AxiomC{}\noLine\UnaryInfC{$[\psi]$}\noLine\UnaryInfC{$\vdots$}\noLine\UnaryInfC{$\alpha$}\RightLabel{$\vee\textsf{E}$} \TrinaryInfC{$\alpha$}\noLine\UnaryInfC{}\DisplayProof\\%\hline
 \multirow{2}{*}{\AxiomC{}\noLine\UnaryInfC{$\phi\vee\psi$} \RightLabel{$\vee\textsf{Com}$}\UnaryInfC{$\psi\vee\phi$}\DisplayProof%\quad\AxiomC{}\noLine\UnaryInfC{$(\phi\vee\psi)\vee\chi$} \RightLabel{$\vee\textsf{Ass}$}\UnaryInfC{$\phi\vee(\psi\vee\chi)$}\DisplayProof
 }
 &\AxiomC{$\phi\vee\psi$}\AxiomC{}\noLine\UnaryInfC{$[\phi]$}\noLine\UnaryInfC{$\vdots$}\noLine\UnaryInfC{$\chi$} \RightLabel{$\vee$\textsf{Mon}} \BinaryInfC{$\chi\vee\psi$}\noLine\UnaryInfC{}\DisplayProof\\
%    \multicolumn{2}{|c|}{%\AxiomC{$\phi\vee\bot$} \RightLabel{$\vee\textsf{E}_\bot$}\UnaryInfC{$\phi$}\DisplayProof\quad
%    \AxiomC{}\noLine\UnaryInfC{$\phi\vee(\psi\vvee\chi)$}\RightLabel{\textsf{Dstr}$\vee\!\vvee$}\UnaryInfC{$(\phi\vee\psi)\vvee(\phi\vee\chi)$} \noLine\UnaryInfC{}\DisplayProof}\\
% &\\
  \hline
 \end{tabular}
 
 \vspace{2pt}

\item[Rule for interaction between the two disjunctions:]\

\vspace{2pt}

\begin{tabular}{|C{0.9\linewidth}|}\hline
%\AxiomC{$\phi\vee\bot$} \RightLabel{$\vee\textsf{E}_\bot$}\UnaryInfC{$\phi$}\DisplayProof\quad
    \AxiomC{}\noLine\UnaryInfC{$\phi\vee(\psi\vvee\chi)$}\RightLabel{\textsf{Dis}\,$\vee\!\vvee$}\UnaryInfC{$(\phi\vee\psi)\vvee(\phi\vee\chi)$} \noLine\UnaryInfC{}\DisplayProof\\
\hline
\end{tabular}

\end{description}

The consequence relation $\Gamma\vdash\phi$ for the  deduction system is defined as usual. In particular,  we write $\phi\dashv\vdash\psi$, if both $\phi\vdash\psi$ and $\psi\vdash\phi$ hold.

%We write $\Gamma\vdash_{\PL(\!\vvee\!)}\phi$ (or simply $\Gamma\vdash\phi$), if $\phi$ is derivable from the set $\Gamma$ of formulas by applying the rules of the system of $\PL(\!\vvee\!)$. We write simply $\phi\vdash\psi$ for $\{\phi\}\vdash\psi$. Two formulas $\phi$ and $\psi$ are said to be {\em provably equivalent}, written as $\phi\dashv\vdash\psi$, if $\phi\vdash\psi$ and $\psi\vdash\phi$.

\end{defn}

%\todo{comment on associativity of $\vee$ -- not needed}
Our rules for negation are more general than those in the (folklore) system in \cite{CiardelliIemhoffYang2020} in that we adopt the standard introduction and elimination rules, while in \cite{CiardelliIemhoffYang2020} the principal formula $\phi$ has to be assumed classical. The soundness of our $\neg$\textsf{I} and $\neg$\textsf{E} rules follows easily from the fact that $\neg\phi$ is flat. The soundness of the other rules are easy to verify.

Since our logic and deduction system are both conservative over classical propositional logic, the above system is clearly complete over classical formulas:

%Restricted to classical formulas, our system coincide with the standard system for classical logic. Thus, 
%Our system is conservative over the standard system of classical propositional logic in the sense of the following lemma:
%The completeness theorem for our system over classical formulas follows easily from the fact that our logic and deduction system are both conservative over classical propositional logic.
\begin{lem}\label{pd_is_cpl_4_classical_fm}
For any set $\Delta\cup\{\alpha\}$ of classical formulas,
\[\Delta\models\alpha\iff  \Delta\vdash\alpha.\]
\end{lem}
\begin{proof}
By flatness of classical formulas and by the equivalence (\ref{flat_singleton_singleval}), it can be verified $\Delta\models\alpha\iff \Delta\models^c\alpha$, where $\models^c$ denotes the semantic entailment relation for the usual (single valuation-based) classical logic \PL. By the completeness theorem of classical propositional logic, we have $\Delta\models^c\alpha\iff \Delta\vdash_{\PL}\alpha$. Finally, since when restricted to classical formulas, our system coincide with the standard system for classical logic, we have $\Delta\vdash_{\PL}\alpha\iff \Delta\vdash_{\PL(\!\vvee\!)}\alpha$.
%
%By Corollary \ref{CPL_entailment_eq_team_entailment} and the completeness theorem of classical propositional logic, we have $\Delta\models\alpha\iff \Delta\models^c\alpha\iff \Delta\vdash_{\PL}\alpha\iff\Delta\vdash_{\PL(\!\vvee\!)}\alpha$, where $\vdash_{\PL}$ is the entailment relation for the standard system of \PL.
\qed\end{proof}

%Thus, all classical validities and rules derivable in the system, such as $\alpha\vee\neg\alpha$.

Let $\phi[\chi/p]$ denote the formula obtained from $\phi$ by replacing one particular instance of $p$ with $\chi$.
An easy induction shows that, as usual, the replacement lemma holds:%\todo{name? replacement?}
%\begin{lem}\label{repl_syn}
\begin{equation}\label{repl_syn}
\text{If }\chi\dashv\vdash\chi'\text{, then }\phi[\chi/p]\vdash\phi[\chi'/p],
\end{equation}
%\end{lem}
where in the negation case one may apply 
the following antitone rule that is clear derivable from $\neg$\textsf{E} and $\neg$\textsf{I}:
\[\AxiomC{$\neg \phi$}\AxiomC{}\noLine\UnaryInfC{$[\psi]$}\noLine\UnaryInfC{$\vdots$}\noLine\UnaryInfC{$\phi$} \RightLabel{$\neg\textsf{Anton}$} \BinaryInfC{$\neg\psi$}\noLine\UnaryInfC{}\DisplayProof\]

%By a standard argument using $\vee\textsf{E}$ and $\vee\textsf{I}$, one can show that the converse direction of the rule \textsf{Dstr}$\vee\!\vvee$ is derivable. Thus, $\phi\vee(\psi\vvee\chi)\dashv\vdash (\phi\vee\psi)\vvee(\phi\vee\chi)$.

As in \cite{CiardelliIemhoffYang2020},  every formula is provably equivalent in the system of $\PL(\vvee)$ to one in  disjunctive normal form (c.f. the more detailed disjunctive normal form  in \Cref{sec:dw-exp}):

\begin{lem}\label{dnf_derivable_plvvee}
For any formula $\phi$, we have that
\[\phi\dashv\vdash\bigvvee_{i\in I}\alpha_i\]
for some finite set $\{\alpha_i\mid i\in I\}$ of classical formulas.
\end{lem}
\begin{proof}
The lemma is proved by induction on $\phi$. The proof for all cases except for negation case can be found  in, e.g., \cite{Ciardelli2015,CiardelliIemhoffYang2020}. Notably, for the key case $\phi=\psi\vee\chi$, one applies the rule \textsf{Dis}$\vee\!\vvee$ along with the induction hypothesis. 

For the new case $\phi=\neg \psi$, by induction hypothesis,  $\psi\dashv\vdash\bigvvee_{i\in I}\alpha_i$ for some classical formulas $\alpha_i$. Thus, 
%$\neg\textsf{Anton}$, 
we have $\neg\psi\dashv\vdash \neg \bigvvee_{i\in I}\alpha_i$ (by (\ref{repl_syn})). We now show $ \neg \bigvvee_{i\in I}\alpha_i\dashv\vdash\bigwedge_{i\in I}\neg\alpha_i$, where the formula on the right hand side is classical.

For the left to right direction, for each $j\in I$, we have 
$\alpha_j\vdash \bigvvee_{i\in I}\alpha_i$ by $\vvee$\textsf{I}. Thus, $\neg \bigvvee_{i\in I}\alpha_i,\alpha_j\vdash \bot$, which implies $\neg \bigvvee_{i\in I}\alpha_i\vdash\neg\alpha_j$ for each $j\in I$ by $\neg$\textsf{I}. Thus,  we conclude $\neg \bigvvee_{i\in I}\alpha_i\vdash\bigwedge_{i\in I}\neg\alpha_i$ by $\wedge$\textsf{I}.

For the converse direction, we first derive $\bigvvee_{i\in I}\alpha_i,\bigwedge_{i\in I}\neg \alpha_i\vdash\bot$ by applying $\vvee$\textsf{E} and $\wedge$\textsf{E}. Then we conclude $\bigwedge_{i\in I}\neg\alpha_i\vdash \neg \bigvvee_{i\in I}\alpha_i$ by $\neg$\textsf{I}.
\qed\end{proof}

The completeness theorem of our system can be proved through an argument that uses the disjunctive normal form result above:

\begin{thm}[Completeness]\label{plvvee_completeness_thm}
For any set  $\Gamma\cup \{\phi\}$ of formulas of $\PL(\vvee)$,
\[\Gamma \models \phi \iff      \Gamma\vdash  \phi .\]
\end{thm} 
\begin{proof}
The detailed proof can be found  in, e.g., \cite{Ciardelli2015,CiardelliIemhoffYang2020}; we only give a proof sketch. Since  $\PL(\vvee)$ is compact (see \cite{VY_PD}), we may assume $\Gamma$ is finite and $\psi=\bigwedge\Gamma$. Suppose $\psi\models \phi$. W.l.o.g., assume further that $\psi=\bigvvee_{i\in I}\alpha_i$ and $\phi=\bigvvee_{j\in J}\beta_j$ are in disjunctive normal form. Then for each $i\in I$, there exists $j_i\in J$ such that $\alpha_i\models\beta_{j_i}$. By \Cref{pd_is_cpl_4_classical_fm}, it follows that $\alpha_i\vdash\beta_{j_i}$. It then follows easily that $\bigvvee_{i\in I}\alpha_i\vdash\bigvvee_{j\in J}\beta_j$, which then gives $\psi\vdash\phi$ by  \Cref{dnf_derivable_plvvee}.
\qed\end{proof}

%\todo{comment that \AxiomC{}\noLine\UnaryInfC{$(\phi\vee\psi)\vee\chi$} \RightLabel{$\vee\textsf{Ass}$}\UnaryInfC{$\phi\vee(\psi\vee\chi)$}\DisplayProof is not needed, and many other rules are simplified}\todo{excluded middle is derivable}

As we demonstrated above, %including the full (intutionistic) negation in the logic $\PL(\vvee)$ does not increase the complexity of its proof theory. B
both the deduction system of $\PL(\vvee)$ with full (intuitionistic) negation and its completeness proof can be easily adapted from those for $\PL(\vvee)$ with restricted negation in the literature. We now summarise a few further folklore observations regarding this system of $\PL(\vvee)$.

%make some further remarks about this folklore system of $\PL(\vvee)$.

First, in the \textsf{RAA} and $\vee$\textsf{E} rules, we restrict the formula $\alpha$ to be classical, but the two rules are sound as long as the formula $\alpha$ is flat. A larger syntactically defined class of flat formula is the class of {\em Harrop formulas}, defined as
\[\delta::=p\mid \bot\mid \neg\phi\mid\delta\wedge\delta\mid \delta\vee\delta\]
where $\phi$ ranges over arbitrary $\PL(\vvee)$ formula. For instance, all classical formulas, $\neg (p\vvee \neg p)$, and $\neg (\phi\vvee\psi)\vee (q\wedge r)$ are  Harrop formulas, whereas $p\vvee \neg p$ is not a Harrop formula.  Harrop formulas are studied originally in the literature on intuitionistic logic \cite{Harrop56}. They are considered by Pun\u{c}och\'{a}\u{r} \cite{FergusonPuncochar2013,Puncochar2016} also in the context of inquisitive logic and team-based logics. (Both contexts consider Harrop formulas in languages with intuitionistic implication, which is not available in $\PL(\vvee)$.) %Harrop formulas are of particular interest in our context as well, as they are all easily verified (by induction) to be flat, even though they can contain the global disjunction $\vvee$ (and are thus not classical).
It is easy to verify by induction that all Harrop formulas are flat, and  by \Cref{pl-exp-comp}, every Harrop formula $\delta$ is equivalent to some classical formula $\delta^\ast$. In fact, $\delta^\ast$ can be obtained directly from $\delta$ by replacing every global disjunction $\vvee$ in $\delta$ with the local disjunction $\vee$ (follows from an easy inductive proof by using the flatness of $\delta$). The more general version of the (sound) rules \textsf{RAA} and $\vee$\textsf{E} with $\alpha$ being a Harrop formula is, obviously, derivable in our (complete) system.

Next, %we give some comments on the earlier presentations of the system above. %, and also some interesting applications of the system.
we remark that some earlier presentations of the system (such as \cite{Ciardelli2015,VY_PD}) include more rules for the local disjunction, particularly the associativity rule  %for the local disjunction:
\[\AxiomC{}\noLine\UnaryInfC{$(\phi\vee\psi)\vee\chi$} \RightLabel{$\vee\textsf{Ass}$}\UnaryInfC{$\phi\vee(\psi\vee\chi)$}\DisplayProof \]
and the rule $\phi\vee\bot/\phi$. 
These rules are actually derivable in the system, via the disjunctive normal form lemma above: %Below we provide this derivation. %We also derive another unnecessary rule $\phi\vee\bot/\phi$ that was included in, e.g., the earlier system in \cite{VY_PD}.

\begin{prop}\label{der_rules}
The following are derivable:
\begin{enumerate}[(i)]
\item $(\phi\vee\psi)\vee\chi\vdash\phi\vee(\psi\vee\chi)$
\item $\phi\vee\bot\vdash\phi$
\end{enumerate}
\end{prop}
\begin{proof}
(i) By \Cref{dnf_derivable_plvvee}, we have that $\phi\dashv\vdash\bigvvee_{i\in I}\alpha_i$, $\psi\dashv\vdash\bigvvee_{j\in J}\beta_j$ and $\chi\dashv\vdash\bigvvee_{k\in K}\gamma_k$ for some classical formulas $\alpha_i$, $\beta_j$ and $\gamma_k$. Then we have that%Then by applying $\vee$\textsf{Mon}, Dis$\vee\vvee$ and classical rules, we have\todo{double check, not quite correct}
\begin{align*}
(\phi\vee\psi)\vee\chi&\vdash\big((\bigvvee_{i\in I}\alpha_i)\vee\bigvvee_{j\in J}\beta_j\big)\vee\bigvvee_{k\in K}\gamma_k\vdash\big( \bigvvee_{i\in I,j\in J}(\alpha_i\vee\beta_j)\big)\vee\bigvvee_{k\in K}\gamma_k\tag{by $\vee$\textsf{Mon}, Dis$\vee\vvee$}\\
&\vdash \bigvvee_{i\in I,j\in J,k\in K}\big((\alpha_i\vee\beta_j)\vee\gamma_k)\vdash \bigvvee_{i\in I,j\in J,k\in K}(\alpha_i\vee(\beta_j\vee\gamma_k)) \tag{by \Cref{pd_is_cpl_4_classical_fm}, classical rule}%\tag{\Cref{pd_is_cpl_4_classical_fm}}
\\
&\vdash (\bigvvee_{i\in I}\alpha_i)\vee\bigvvee_{j\in J,k\in K}(\beta_j\vee\gamma_k)\vdash  (\bigvvee_{i\in I}\alpha_i)\vee\big((\bigvvee_{j\in J}\beta_j)\vee\bigvvee_{k\in K}\gamma_k\big)\tag{by $\vee$\textsf{Mon}, Dis$\vee\vvee$}\\
&\vdash\phi\vee(\psi\vee\chi).
\end{align*}

(ii) By \Cref{dnf_derivable_plvvee}, we have that $\phi\dashv\vdash\bigvvee_{i\in I}\alpha_i$ for some classical formulas $\alpha_i$. Then,%Now,  we have 
\begin{align*}
\phi\vee\bot&\vdash (\bigvvee_{i\in I}\alpha_i)\vee\bot\tag{$\vee$\textsf{Mon}}\\
&\vdash\bigvvee_{i\in I}(\alpha_i\vee\bot)\tag{\textsf{$\vee$\textsf{Com}, Dis$\vee\vvee$}}\\
&\vdash \bigvvee_{i\in I}\alpha_i\tag{By \Cref{pd_is_cpl_4_classical_fm}, classical rule}\\
&\vdash\phi.
\end{align*}
\qed\end{proof}

These derived rules allow us to derive  interesting clauses such as the following %(classically also sound) 
ones:
\begin{prop}\label{deri_rule}
The following are derivable:
\begin{enumerate}[(i)]
\item $\alpha\vee\phi,\neg\alpha\vdash\phi$
\item $\alpha\vee\phi,\neg\alpha\vee\psi\vdash \phi\vee\psi$
\end{enumerate}
\end{prop}
\begin{proof}
(i) Since $\alpha,\neg\alpha\vdash\bot$ (by $\neg$\textsf{E}), we have $\alpha\vee\phi,\neg\alpha\vdash \bot\vee\phi$ by $\vee$\textsf{Mon}. By \Cref{der_rules}(i), we have $\bot\vee\phi\vdash\phi$. Putting these together, we conclude $\alpha\vee\phi,\neg\alpha\vdash\phi$. %as desired.

(ii) Follows from (i), and $\vee$\textsf{Mon} applied to $\neg\alpha$ in the disjunction $\neg\alpha\vee\psi$.
\qed\end{proof}

%Lastly, let us point out that excluded middle law for classical formulas $\vdash\neg\alpha\vee\alpha$ follows from \Cref{pd_is_cpl_4_classical_fm}. Moreover, since our system contains the standard introduction and elimination rule for negation (i.e., $\neg$\textsf{I} and $\neg$\textsf{E}), and \textsf{RAA} rule holds for Harrop formulas, a standard derivation also gives us excluded middle law also for Harrop formulas, i.e., $\vdash\neg\delta\vee\delta$.
%
%\todo{extensions}

Let us end with some comments on the system for %variants and extensions of the system for other downward closed languages. 
%One may also consider 
the language of $\PL(\vvee)$ extended with the intuitionistic implication $\to$, which can also be viewed as inquisitive logic extended with the local disjunction. The (classical) systems in  \cite{Ciardelli2015,CiardelliIemhoffYang2020} for this language contain two additional rules for intuitionistic implication and the \textsf{Split} rule:

\begin{description}
\item[Rules for intuitionistic implication:]\

\vspace{2pt}

{\normalfont %\begin{center}
\renewcommand{\arraystretch}{1.8}
\begin{tabular}{|C{0.55\linewidth}C{0.55\linewidth}|}
\hline
\AxiomC{}\noLine \UnaryInfC{[$\phi$]}\noLine \UnaryInfC{$\vdots$}\noLine\UnaryInfC{$\psi$}\RightLabel{$\to\textsf{I}$}\UnaryInfC{$\phi\to\psi$}\noLine\UnaryInfC{}\DisplayProof
&

\AxiomC{$\phi\to \psi$} \AxiomC{$\phi$}\RightLabel{$\to\mathsf{E}$}\BinaryInfC{$\psi$}\noLine\UnaryInfC{}\DisplayProof\\
\hline
\end{tabular}
%\end{center}
}

\vspace{2pt}

\item[Split rule:]\

\vspace{2pt}

{\normalfont %\begin{center}
\renewcommand{\arraystretch}{1.8}
\begin{tabular}{|C{\linewidth}|}
\hline
\AxiomC{}\noLine \UnaryInfC{$\alpha\to(\phi\vvee\psi)$}\RightLabel{\textsf{Split}}\UnaryInfC{$(\alpha\to\phi)\vvee(\alpha\to\psi)$}\noLine\UnaryInfC{}\DisplayProof
\\
\hline
\end{tabular}
%\end{center}
}
\end{description}

%Split axiom/rule is derivable from $\vvee\vee$Dis: $\alpha\to(\phi\vvee\psi)/(\alpha\to\phi)\vvee(\alpha\to\psi)$

Intuitionistic implications of the form $\alpha\to\phi$ with $\alpha$ classical behaves classically in the sense that it is equivalent to the local disjunction $\neg\alpha\vee\phi$. This fact can also be obtained proof-theoretically in our extended deduction system:

\begin{prop}\label{imp_defi}
$\alpha\to\phi\dashv\vdash\neg\alpha\vee\phi$.%\todo{holds for $\alpha$ as well}
\end{prop}
\begin{proof}
For the left to right direction, note that we have %that 
$\vdash\neg\alpha\vee\alpha$ by \Cref{pd_is_cpl_4_classical_fm}. Now, $\neg\alpha\vee\alpha,\alpha\to\phi\vdash \neg\alpha\vee\phi$ by $\to$\textsf{E} and $\vee$\textsf{Mon}. Putting these together, we obtain $\alpha\to\phi\vdash \neg\alpha\vee\phi$.

For the right to left direction, by \Cref{deri_rule}(i) we have $\neg\alpha\vee\phi,\neg\neg\alpha\vdash \phi$, and by  \Cref{pd_is_cpl_4_classical_fm}  we have $\alpha\vdash\neg\neg\alpha$. Thus, we obtain $\neg\alpha\vee\phi,\alpha\vdash \phi$, which by $\to$\textsf{I} implies $\neg\alpha\vee\phi\vdash\alpha\to\phi$.
\qed\end{proof}

Given this, the \textsf{Split} rule, which is often viewed as a distinguishing feature of inquisitive logic, is in fact derivable from \textsf{Dis}$\vvee\vee$ and other rules in the  system. %We illustrate this fact with a more general version of the \textsf{Split} rule $\delta\to(\phi\vvee\psi)/(\delta\to \phi)\vvee(\delta\to \psi)$ with $\delta$ being a Harrop formula. One can easily verify that the rule is sound. 
To see how to derive it in the system, 
%To see why, 
note that 
%The \textsf{Split} rule can be viewed as an distinguishing feature of inquisitive logic.
%
%\textsf{Split} rule is in fact derivable from \textsf{Dis}$\vvee\vee$ and other rules in the  system.
%\begin{prop}
%$\alpha\to\phi\vvee\psi\vdash(\alpha\to\phi)\vvee(\alpha\to\psi)$.
%\end{prop}
%\begin{proof}
by \Cref{imp_defi}, to derive the \textsf{Split} rule it is sufficient to derive
\[\neg\alpha\vee(\phi\vvee\psi)\vdash (\neg\alpha\vee\phi)\vvee(\neg\alpha\vee\psi),\] but this  follows immediately from \textsf{Dis}$\vee\vvee$. %Clearly, this derivation also applies to a more general version of the \textsf{Split} rule $\delta\to(\phi\vvee\psi)/(\delta\to \phi)\vvee(\delta\to \psi)$ with $\delta$ being a Harrop formula.\todo{cite Vit for this Harrop version of the Split}

\section{Concluding remarks and further directions}

In this paper, we have explored propositional team-based logics with full (intuitionistic) negation. We illustrated how the complement of certain properties can be characterised already by other means without directly using negatoin, % than simply negating the properties (with classical negation), 
and we also showed that including the full intuitionistic negation does not substantially affect the natural deduction system of the downward closed logic $\PL(\vvee)$. The original version of $\PL(\vvee)$ with negation restricted to classical formulas only has been recently introduced  a labelled sequent calculus \cite{BarberoGirlandoMullerYang2024} and a deep inference style calculus \cite{AnttilaIemhoffYang2024}. Generalizing these calculi to cover full intuitionistic negation as well is left as future work.

Another downward closed logic with full intuitionistic negation we considered in the paper is $\PL(\dep(\cdot))$. The version of $\PL(\dep(\cdot))$ with atomic negation was axiomatized in \cite{VY_PD} with a complete natural deduction system that contains all rules for $\PL(\vvee)$ (with restricted negation) except for those that involve global disjunction $\vvee$ and some additional (deep inference style) rules for dependence atoms. It is essentially clear that a complete natural deduction system for $\PL(\dep(\cdot))$ with full intuitionistic negation can be obtained in the same way as we did in \Cref{plvvee_deduction_system} for $\PL(\vvee)$. We will leave the proof details for future research, as they are more involved. 

Similarly, the modal versions of the logics  we considered in this paper with restricted negation have been axiomatized in \cite{Yang17MD}.  Generalizing these systems to encompass the full negation case should be reasonably achievable. %, even though the interaction between negation and modalities potentially make the details more involved. 
Working out the details, especially the interaction between negation and modalities, is left as future work.

%\todo{modal logic? Everything here can also be generalised to the modal case, though the interaction between negation and modalities should be paid more attention. We leave the details also for future research.}

% from the system of $\PL(\vvee)$ with full intuitionistic negation as presented in \Cref{plvvee_deduction_system}, by using $\neg(\phi\vvee\psi)\equiv \neg\phi\wedge\neg\psi$. We leave the details for further research.

%\todo{for $\PL(\dep(\cdot))$, comment on deep inference style rules, and remark that one can indeed develop deep inference style sequent calculi, see upcoming paper of my paper with Aleksi and Rosalie}\todo{extend also the sequent style calculi to cover full negation case}

Much of the flexibility in these results relies on the downward closed setting. As for the non-downward closed logics we considered, it is not clear whether, e.g., a direct adaptation to the system of $\PL(\subseteq)$ with restricted negation gives rise to a complete system for $\PL(\subseteq)$ with full negation as well. In particular, the completeness proof for the system of $\PL(\vvee)$ relies heavily on a syntactic normal form lemma (\Cref{dnf_derivable_plvvee}), which, as we showed, can be easily generalized to cover the full negation case. However, the corresponding syntactic normal form lemma for  $\PL(\subseteq)$ from \cite{Yang2022} involves more complex formulas, and thus generalizing this to cover the full negation case does not seem to be straightforward.

Finally, we remark that apart from the intuitionstic implication $\to$, the full negation induces two other implications defined in terms of the two different versions of disjunctions:
\[\phi\rightharpoondown\psi:=\neg\phi\vee\psi\quad\text{ and }\quad \phi\barrightharpoon\psi:=\neg\phi\vvee\psi.\]
It is easy to verify that $\phi%\rightharpoonup
\rightharpoondown\psi\equiv\phi\to\psi$ if $\phi$ is union closed. Modus Ponens is sound for both implications, but deduction theorem fails in general for both implications. %Recent work by Barbero \cite{?} showed that there is no implication in such setting that satisfies both Modus Ponens and deduction theorem.
%($\Gamma,\phi\models\psi$ implies $\Gamma\models\phi\rightharpoondown\psi$). 
 Studying the properties of these implications and possibly other connectives induced by the full negation is an interesting further direction.

%In terms of the De Morgan's law, we have the following:
%\[\neg(\phi\vvee\psi)\equiv \neg \phi\wedge \neg\psi\equiv \neg (\phi\vee\psi),\quad \neg (\phi\wedge\psi)\equiv \neg\phi\vee\neg\psi.\]

%\todo{are negation rules sound for non-downward closed logics? If so, checking those systems are further directions}\todo{still sound, but may not be complete. Negation of the normal form does not behave well. How to axiomatize those logic belongs to further directions.}\todo{but this is not surprising as in the non-downward closed setting, already intuitionistic implication does not behave well. It was shown by Fausto that there is no implication in such setting that satisfies deduction theorem.}

\begin{credits}
\subsubsection{\ackname} The author would like to thank Aleksi Anttila and Fausto Barbero for discussions related to the topic of this paper.

%\subsubsection{\discintname}
%It is now necessary to declare any competing interests or to specifically
%state that the authors have no competing interests. Please place the
%statement with a bold run-in heading in small font size beneath the
%(optional) acknowledgments\footnote{If EquinOCS, our proceedings submission
%system, is used, then the disclaimer can be provided directly in the system.},
%for example: The authors have no competing interests to declare that are
%relevant to the content of this article. Or: Author A has received research
%grants from Company W. Author B has received a speaker honorarium from
%Company X and owns stock in Company Y. Author C is a member of committee Z.
\end{credits}

\bibliographystyle{splncs04}

%\bibliography{../fan}

\begin{thebibliography}{10}
\providecommand{\url}[1]{\texttt{#1}}
\providecommand{\urlprefix}{URL }
\providecommand{\doi}[1]{https://doi.org/#1}

\bibitem{AbVan09}
Abramsky, S., V\"{a}\"{a}n\"{a}nen, J.: From {IF} to {BI}. Synthese
  \textbf{167}(2),  207--230 (2009)

\bibitem{Anttila2024}
Anttila, A.: Further Remarks on the Dual Negation in Team Logics (2024)

\bibitem{AnttilaIemhoffYang2024}
Anttila, A., Iemhoff, R., Yang, F.: Deep-inference sequent calculi for
  propositional team-based logics (2024)

\bibitem{BarberoGirlandoMullerYang2024}
Barbero, F., Girlando, M., M\"{u}ller, V., Yang, F.: Labelled Proof Systems for
  Inquisitive and Team-based logics (2024)

\bibitem{Burgess_negation_03}
Burgess, J.P.: A remark on henkin sentences and their contraries. Notre Dame
  Journal of Formal Logic  \textbf{44}(3),  185--188 (2003)

\bibitem{Ciardelli2015}
Ciardelli, I.: Dependency as question entailment. In: Abramsky, S., Kontinen,
  J., V{\"a}{\"a}n{\"a}nen, J., Vollmer, H. (eds.) Dependence Logic: Theory and
  Application, pp. 129--182. Progress in Computer Science and Applied Logic,
  Birkhauser (2016)

\bibitem{CiardelliIemhoffYang2020}
Ciardelli, I., Iemhoff, R., Yang, F.: Questions and dependency in
  intuitionistic logic. Notre Dame Journal of Formal Logic  \textbf{61}(1),
  75--115 (2020)

\bibitem{InquiLog}
Ciardelli, I., Roelofsen, F.: Inquisitive logic. Journal of Philosophical Logic
   \textbf{40}(1),  55--94 (2011)

\bibitem{FergusonPuncochar2013}
Ferguson, T., Pun\u{c}och\'{a}\u{r}, V.: Structural completeness and
  superintuitionistic inquisitive logics. In: Helle Hvid~Hansen,
  A.S..R.J.G.B.D.Q. (ed.) Logic, Language, Information, and Computation: 29th
  International Workshop (WoLLIC). pp. 194--210 (2023)

\bibitem{Pietro_I/E}
Galliani, P.: Inclusion and exclusion in team semantics: On some logics of
  imperfect information. Annals of Pure and Applied Logic  \textbf{163}(1),
  68--84 (January 2012)

\bibitem{inclusion_logic_GH}
Galliani, P., Hella, L.: Inclusion logic and fixed point logic. In: Computer
  Science Logic 2013. Leibniz International Proceedings in Informatics
  (LIPIcs), vol.~23, pp. 281--295. Schloss Dagstuhl - Leibniz-Zentrum fuer
  Informatik (2013)

\bibitem{D_Ind_GV}
Gr{\"a}del, E., V{\"a}{\"a}n{\"a}nen, J.: Dependence and independence. Studia
  Logica  \textbf{101}(2),  399--410 (April 2013)

\bibitem{Hannula_fo_ind_13}
Hannula, M.: Axiomatizing first-order consequences in independence logic.
  Annals of Pure and Applied Logic  \textbf{166}(1),  61--91 (2015)

\bibitem{Harrop56}
Harrop, R.: On disjunctions and existential statements in intuitionistic
  systems of logic. Mathematische Annalen  \textbf{132}(4),  347--361 (1956)

\bibitem{HLSV14}
Hella, L., Luosto, K., Sano, K., Virtema, J.: The expressive power of modal
  dependence logic. In: Proceedings of Advances in Modal Logic 10. pp.
  294--312. College Publications (2014)

\bibitem{HellaStumpf15}
Hella, L., Stumpf, J.: The expressive power of modal logic with inclusion
  atoms. In: Proceedings of the 6th GandALF. pp. 129--143 (2015)

\bibitem{henkin61}
Henkin, L.: Some remarks on infinitely long formulas. In: Infinitistic Methods.
  pp. 167--183. Proceedings Symposium Foundations of Mathematics, Pergamon,
  Warsaw (1961)

\bibitem{Hintikka98book}
Hintikka, J.: The Principles of Mathematics Revisited. Cambridge University
  Press (1998)

\bibitem{HinSan96}
Hintikka, J., Sandu, G.: Game-theoretical semantics. In: van Benthem, J., ter
  Meulen, A. (eds.) Handbook of Logic and Language. Elsevier (1996)

\bibitem{Hodges1997a}
Hodges, W.: Compositional semantics for a language of imperfect information.
  Logic Journal of the IGPL  \textbf{5},  539--563 (1997)

\bibitem{Hodges1997b}
Hodges, W.: Some strange quantifiers. In: Mycielski, J., Rozenberg, G.,
  Salomaa, A. (eds.) Structures in Logic and Computer Science: A Selection of
  Essays in Honor of A. Ehrenfeucht, Lecture Notes in Computer Science,
  vol.~1261, pp. 51--65. London: Springer (1997)

\bibitem{KontinenNurmi2011}
Kontinen, J., Nurmi, V.: Team logic and second-order logic. Fundamenta
  Informaticae  \textbf{106},  259--272 (2011)

\bibitem{KontVan09}
Kontinen, J., V\"{a}\"{a}n\"{a}nen, J.: On definability in dependence logic.
  Journal of Logic, Language and Information  \textbf{18(3)},  317--332
  (Erratum: \emph{the same journal, 20(1)} (2011), 133--134) (2009)

\bibitem{Negation_D_KV}
Kontinen, J., V{\"a}{\"a}n{\"a}nen, J.: A remark on negation in dependence
  logic. Notre Dame Journal of Formal Logic  \textbf{52}(1),  55--65 (2011)

\bibitem{Axiom_fo_d_KV}
Kontinen, J., V{\"a}{\"a}n{\"a}nen, J.: Axiomatizing first-order consequences
  in dependence logic. Annals of Pure and Applied Logic  \textbf{164}(11)
  (2013)

\bibitem{Luck18}
L\"{u}ck, M.: Axiomatizations of team logics. Annals of Pure and Applied Logic
  \textbf{169}(9),  928--969 (2018)

\bibitem{Puncochar15}
Pun\u{c}och\'{a}\u{r}, V.: Weak negation in inquisitive semantics. Journal of
  Logic, Language and Information  \textbf{24},  323--355 (2015)

\bibitem{Puncochar2016}
Pun\u{c}och\'{a}\u{r}, V.: Journal of philosophical logic. A Generalization of
  Inquisitive Semantics  \textbf{45}(4),  399--428 (August 2016)

\bibitem{Van07dl}
V{\"a}{\"a}n{\"a}nen, J.: Dependence Logic: A New Approach to Independence
  Friendly Logic. Cambridge: Cambridge University Press (2007)

\bibitem{Yang2011}
Yang, F.: Expressing second-order sentences in intuitionistic dependence logic.
  Studia Logica  \textbf{101}(2),  323--342 (2013)

\bibitem{Yang17MD}
Yang, F.: Modal dependence logics: Axiomatizations and model-theoretic
  properties. Logic Journal of the IGPL  \textbf{25}(5),  773--805 (October
  2017)

\bibitem{Yang15}
Yang, F.: Uniform definability in propositional dependence logic. Review of
  Symbolic Logic  (2017)

\bibitem{Yang_neg18}
Yang, F.: Negation and partial axiomatizations of dependence and independence
  logic revisited. Annals of Pure and Applied Logic  \textbf{170}(9),
  1128--1149 (September 2019)

\bibitem{YangInc20}
Yang, F.: Axiomatizing first-order consequences in inclusion logic.
  Mathematical Logic Quarterly  \textbf{66}(2),  195--216 (July 2020)

\bibitem{Yang2022}
Yang, F.: Propositional union closed team logics. Annals of Pure and Applied
  Logic  \textbf{173}(6),  103102 (June 2022)

\bibitem{VY_PD}
Yang, F., V\"{a}\"{a}n\"{a}nen, J.: Propositional logics of dependence. Annals
  of Pure and Applied Logic  \textbf{167}(7),  557--589 (July 2016)

\bibitem{YangVaananen:17PT}
Yang, F., V{\"a}{\"a}n{\"a}nen, J.: Propositional team logics. Annals of Pure
  and Applied Logic  \textbf{168}(7),  1406--1441 (July 2017)

\end{thebibliography}
%Author, F.: Article title. Journal \textbf{2}(5), 99--110 (2016)
%Author, F., Author, S.: Title of a proceedings paper. In: Editor,
%F., Editor, S. (eds.) CONFERENCE 2016, LNCS, vol. 9999, pp. 1--13.
%Springer, Heidelberg (2016). \doi{10.10007/1234567890}
%Author, F., Author, S., Author, T.: Book title. 2nd edn. Publisher,
%Author, A.-B.: Contribution title. In: 9th International Proceedings
%on Proceedings, pp. 1--2. Publisher, Location (2010)
%LNCS Homepage, \url{http://www.springer.com/lncs}, last accessed 2023/10/25
%\end{thebibliography}
\end{document}